# MOMENT INEQUALITIES FOR FUNCTIONS OF INDEPENDENT RANDOM VARIABLES

By Stéphane Boucheron[1], Olivier Bousquet, Gábor Lugosi[2]
and Pascal Massart

*CNRS-Université Paris-Sud, Max Planck Institute, Pompeu Fabra
University and Université Paris-Sud*

A general method for obtaining moment inequalities for functions of independent random variables is presented. It is a generalization of the entropy method which has been used to derive concentration inequalities for such functions [Boucheron, Lugosi and Massart *Ann. Probab.* **31** (2003) 1583–1614], and is based on a generalized tensorization inequality due to Latała and Oleszkiewicz [*Lecture Notes in Math.* **1745** (2000) 147–168]. The new inequalities prove to be a versatile tool in a wide range of applications. We illustrate the power of the method by showing how it can be used to effortlessly re-derive classical inequalities including Rosenthal and Kahane–Khinchine-type inequalities for sums of independent random variables, moment inequalities for suprema of empirical processes and moment inequalities for Rademacher chaos and $U$-statistics. Some of these corollaries are apparently new. In particular, we generalize Talagrand's exponential inequality for Rademacher chaos of order 2 to any order. We also discuss applications for other complex functions of independent random variables, such as suprema of Boolean polynomials which include, as special cases, subgraph counting problems in random graphs.

**1. Introduction.** During the last twenty years, the search for upper bounds for exponential moments of functions of independent random variables, that is, for concentration inequalities, has been a flourishing area of probability theory. Recent developments in random combinatorics, statistics and empir-

Received February 2003; revised January 2004.

[1]Supported by EU Working Group RAND-APX, binational PROCOPE Grant 05923XL.

[2]Supported by the Spanish Ministry of Science and Technology and FEDER, Grant BMF2003-03324.

ical process theory have prompted the search to moment inequalities dealing with possibly nonexponentially integrable random variables.

Paraphrasing Talagrand in [41], we may argue that

> While Rosenthal–Pinelis inequalities for higher moments of sums of independent random variables are at the core of classical probabilities, there is a need for new abstract inequalities for higher moments of more general functions of many independent random variables.

The aim of this paper is to provide such general-purpose inequalities. Our approach is based on a generalization of Ledoux's entropy method (see [26, 28]). Ledoux's method relies on abstract functional inequalities known as logarithmic Sobolev inequalities and provides a powerful tool for deriving exponential inequalities for functions of independent random variables; see [6, 7, 8, 14, 30, 31, 36] for various applications. To derive moment inequalities for general functions of independent random variables, we elaborate on the pioneering work of Latała and Oleszkiewicz [25] and describe so-called $\phi$-Sobolev inequalities which interpolate between Poincaré's inequality and logarithmic Sobolev inequalities (see also [4] and Bobkov's arguments in [26]).

This paper proposes general-purpose inequalities for polynomial moments of functions of independent variables. Many of the results parallel those obtained in [7] for exponential moments, based on the entropy method. In fact, the exponential inequalities of [7] may be obtained (up to constants) as corollaries of the results presented here.

Even though the new inequalities are designed to handle very general functions of independent random variables, they prove to be surprisingly powerful in bounding moments of well-understood functions such as sums of independent random variables and suprema of empirical processes. In particular, we show how to apply the new results to effortlessly re-derive Rosenthal and Kahane–Khinchine-type inequalities for sums of independent random variables, Pinelis' moment inequalities for suprema of empirical processes and moment inequalities for Rademacher chaos. Some of these corollaries are apparently new. Here we mention Theorem 14 which generalizes Talagrand's (upper) tail bound [40] for Rademacher chaos of order 2 to Rademacher chaos of any order. We also provide some other examples such as suprema of Boolean polynomials which include, as special cases, subgraph counting problems in random graphs.

The paper is organized as follows. In Section 2, we state the main results of this paper, Theorems 2–4, as well as a number of corollaries. The proofs of the main results are given in Sections 4 and 5. In Section 4, abstract $\phi$-Sobolev inequalities which generalize logarithmic Sobolev inequalities are introduced. These inequalities are based on a "tensorization property" of certain functionals called $\phi$-entropies. The tensorization property is based on a duality formula, stated in Lemma 1. In Appendix A.1, some further facts are gathered about the tensorization property of $\phi$-entropies.

In Section 6, the main theorems are applied to sums of independent random variables. This leads quite easily to suitable versions of Marcinkiewicz's, Rosenthal's and Pinelis' inequalities. In Section 7, Theorems 2 and 3 are applied to suprema of empirical processes indexed by possibly nonbounded functions, leading to a version of an inequality due to Giné, Latała and Zinn [16] with explicit and reasonable constants. In Section 8, we derive moment inequalities for conditional Rademacher averages. In Section 9, a new general moment inequality is obtained for Rademacher chaos of any order, which generalizes Talagrand's inequality for Rademacher chaos of order 2. We also give a simple proof of Bonami's inequality.

In Section 10, we consider suprema of Boolean polynomials. Such problems arise, for example, in random graph theory where an important special case is the number of small subgraphs in a random graph.

Some of the routine proofs are gathered in the Appendix.

## 2. Main results.

2.1. *Notation.* We begin by introducing some notation used throughout the paper. Let $X_1, \ldots, X_n$ denote independent random variables taking values in some measurable set $\mathcal{X}$. Denote by $X_1^n$ the vector of these $n$ random variables. Let $F : \mathcal{X}^n \to \mathbb{R}$ be some measurable function. We are concerned with moment inequalities for the random variable

$$Z = F(X_1, \ldots, X_n).$$

Throughout, $\mathbb{E}[Z]$ denotes expectation of $Z$ and $\mathbb{E}[Z|\mathcal{F}]$ denotes conditional expectation with respect to $\mathcal{F}$. $X_1', \ldots, X_n'$ denote independent copies of $X_1, \ldots, X_n$, and we write

$$Z_i' = F(X_1, \ldots, X_{i-1}, X_i', X_{i+1}, \ldots, X_n).$$

Define the random variables $V^+$ and $V^-$ by

$$V^+ = \mathbb{E}\left[\sum_{i=1}^n (Z - Z_i')_+^2 | X_1^n\right]$$

and

$$V^- = \mathbb{E}\left[\sum_{i=1}^n (Z - Z_i')_-^2 | X_1^n\right],$$

where $x_+ = \max(x, 0)$ and $x_- = \max(-x, 0)$ denote the positive and negative parts of a real number $x$. The variables $V^+$ and $V^-$ play a central role in [7]. In particular, it is shown in [7] that the moment generating function of $Z - \mathbb{E}Z$ may be bounded in terms of the moment generating functions of $V^+$ and $V_-$. The main results of the present paper relate the moments of $Z$ to lower-order moments of these variables.

In the sequel, $Z_i$ will denote an arbitrary measurable function $F_i$ of $X^{(i)} = X_1, \ldots, X_{i-1}, X_{i+1}, \ldots, X_n$, that is,

$$Z_i = F_i(X_1, \ldots, X_{i-1}, X_{i+1}, \ldots, X_n).$$

Finally, define

$$V = \sum_{i=1}^n (Z - Z_i)^2.$$

Throughout the paper, the notation $\|Z\|_q$ is used for

$$\|Z\|_q = (\mathbb{E}[|Z|^q])^{1/q},$$

where $q$ is a positive number.

Next we introduce two constants used frequently in the paper. Let

$$\kappa = \frac{\sqrt{e}}{2(\sqrt{e} - 1)} < 1.271.$$

Let $\kappa_1 = 1$ and for any integer $q \geq 2$, define

$$\kappa_q = \frac{1}{2}\left(1 - \left(1 - \frac{1}{q}\right)^{q/2}\right)^{-1}.$$

Then $(\kappa_q)$ increases to $\kappa$ as $q$ goes to infinity. Also, define

$$K = \frac{1}{e - \sqrt{e}} < 0.935.$$

2.2. *Basic theorems.* Recall first one of the first general moment inequalities, proved by Efron and Stein [15], and further improved by Steele [37]:

PROPOSITION 1 (Efron–Stein inequality).

$$\mathrm{Var}[Z] \leq \tfrac{1}{2}\mathbb{E}\left[\sum_{i=1}^n (Z - Z_i')^2\right].$$

Note that this inequality becomes an equality if $F$ is the sum of its arguments. Generalizations of the Efron–Stein inequality to higher moments of sums of independent random variables have been known in the literature as Marcinkiewicz's inequalities (see, e.g., [13], page 34). Our purpose is to describe conditions under which versions of Marcinkiewicz's inequalities hold for general functions $F$.

In [7], inequalities for exponential moments of $Z$ are derived in terms of the behavior of $V^+$ and $V^-$. This is quite convenient when exponential moments of $Z$ scale nicely with $n$. In many situations of interest this is not the case, and bounds on exponential moments of roots of $Z$ rather than bounds on exponential moments of $Z$ itself are obtained (e.g., the triangle counting problem in [7]). In such situations, relating the polynomial moments of $Z$ to $V^+$, $V^-$ or $V$ may prove more convenient.

In the simplest settings, $V^+$ and $V^-$ are bounded by a constant. It was shown in [7] that in this case $Z$ exhibits a sub-Gaussian behavior. Specifically, it is shown in [7] that if $V^+ \leq c$ almost surely for some positive constant $c$, then for any $\lambda > 0$,

$$\mathbb{E}e^{\lambda(Z-\mathbb{E}[Z])} \leq e^{\lambda^2 c}.$$

Our first introductory result implies sub-Gaussian bounds for the polynomial moments of $Z$:

THEOREM 1. *If $V^+ \leq c$ for some constant $c \geq 0$, then for all integers $q \geq 2$,*

$$\|(Z - \mathbb{E}[Z])_+\|_q \leq \sqrt{Kqc}.$$

*[Recall that $K = 1/(e - \sqrt{e}) < 0.935$.] If furthermore $V^- \leq c$, then for all integers $q \geq 2$,*

$$\|Z\|_q \leq \mathbb{E}[Z] + 2^{1/q}\sqrt{Kqc}.$$

The main result of this paper is the following inequality.

THEOREM 2. *For any real $q \geq 2$,*

$$\|(Z - \mathbb{E}[Z])_+\|_q \leq \sqrt{\left(1 - \frac{1}{q}\right)2\kappa_q q\|V^+\|_{q/2}}$$
$$\leq \sqrt{2\kappa q\|V^+\|_{q/2}} = \sqrt{2\kappa q}\|\sqrt{V^+}\|_q$$

*and*

$$\|(Z - \mathbb{E}[Z])_-\|_q \leq \sqrt{\left(1 - \frac{1}{q}\right)2\kappa_q q\|V^-\|_{q/2}}$$
$$\leq \sqrt{2\kappa q\|V^-\|_{q/2}} = \sqrt{2\kappa q}\|\sqrt{V^-}\|_q.$$

REMARK. To better understand our goal, recall Burkholder's inequalities [9, 10] from martingale theory. Burkholder's inequalities may be regarded as extensions of Marcinkiewicz's inequalities to sums of martingale increments. They are natural candidates for deriving moment inequalities for a function $Z = F(X_1, \ldots, X_n)$ of many independent random variables. The approach mimics the method of bounded differences (see [32, 33]) classically used to derive Bernstein- or Hoeffding-like inequalities under similar circumstances. The method works as follows: let $\mathcal{F}_i$ denote the $\sigma$-algebra generated by the sequence $(X_1^i)$. Then the sequence $M_i = \mathbb{E}[Z|\mathcal{F}_i]$ is an $\mathcal{F}_i$-adapted martingale (the Doob martingale associated with $Z$). Let $\langle Z \rangle$ denote the associated *quadratic variation*

$$\langle Z \rangle = \sum_{i=1}^{n} (M_i - M_{i-1})^2,$$

let $[Z]$ denote the associated *predictable quadratic variation*

$$[Z] = \sum_{i=1}^{n} \mathbb{E}[(M_i - M_{i-1})^2 | \mathcal{F}_{i-1}],$$

and let $M$ be defined as $\max_{1 \leq i \leq n} |Z_i - Z_{i-1}|$. Burkholder's inequalities [9, 10] (see also [12], page 384) imply that for $q \geq 2$,

$$\|Z - \mathbb{E}[Z]\|_q \leq (q-1)\sqrt{\|\langle Z \rangle\|_{q/2}} = (q-1)\|\sqrt{\langle Z \rangle}\|_q.$$

Note that the dependence on $q$ in this inequality differs from the dependence in Theorem 2. It is known that for general martingales, Burkholder's inequality is essentially unimprovable (see [10], Theorem 3.3). (However, for the special case of Doob martingale associated with $Z$ this bound is perhaps improvable.) The Burkholder–Rosenthal–Pinelis inequality ([34], Theorem 4.1) implies that there exists a universal constant $C$ such that

$$\|Z - \mathbb{E}[Z]\|_q \leq C(\sqrt{q\|[Z]\|_{q/2}} + q\|M\|_q).$$

If one has some extra information on the sensitivity of $Z$ with respect to its arguments, such inequalities may be used to develop a strict analogue of the method of bounded differences (see [33]) for moment inequalities. In principle such an approach should provide tight results, but finding good bounds on the moments of the quadratic variation process often proves quite difficult.

The inequalities introduced in this paper have a form similar to those obtained by resorting to Doob's martingale representation and Burkholder's inequality. But, instead of relying on the quadratic variation process, they rely on a more tractable quantity. Indeed, in many cases $V^+$ and $V^-$ are easier to deal with than $[Z]$ or $\langle Z \rangle$.

Below we present two variants of Theorem 2 which may be more convenient in some applications.

THEOREM 3. *Assume that $Z_i \leq Z$ for all $1 \leq i \leq n$. Then for any real $q \geq 2$,*

$$\|(Z - \mathbb{E}[Z])_+\|_q \leq \sqrt{\kappa_q q \|V\|_{q/2}} \leq \sqrt{\kappa q \|V\|_{q/2}}.$$

Even though Theorem 2 provides some information concerning the growth of moments of $(Z - \mathbb{E}[Z])_-$, this information may be hard to exploit in concrete cases. The following result relates the moments of $(Z - \mathbb{E}[Z])_-$ with $\|V^+\|_q$ rather than with $\|V^-\|_q$. This requires certain boundedness assumptions on the increments of $Z$.

THEOREM 4. *If for some positive random variable $M$,*

$$(Z - Z_i')_+ \leq M \qquad \text{for every } 1 \leq i \leq n,$$

*then for every real $q \geq 2$,*

$$\|(Z - \mathbb{E}[Z])_-\|_q \leq \sqrt{C_1 q(\|V^+\|_{q/2} \vee q\|M\|_q^2)},$$

*where $C_1 < 4.16$. If, on the other hand,*

$$0 \leq Z - Z_i \leq M \qquad \text{for every } 1 \leq i \leq n,$$

*then*

$$\|(Z - \mathbb{E}[Z])_-\|_q \leq \sqrt{C_2 q(\|V\|_{q/2} \vee q\|M\|_q^2)},$$

*where $C_2 < 2.42$.*

2.3. *Corollaries.* Next we derive some general corollaries of the main theorems which provide explicit estimates under various typical conditions on the behavior of $V^+$, $V^-$ or $V$.

The first corollary, obtained from Theorem 3, is concerned with functionals $Z$ satisfying $V \leq Z$. Such functionals were at the center of attention in [6] and [36] where they were called self-bounded functionals. They encompass sums of bounded nonnegative random variables, suprema of nonnegative empirical processes, configuration functions in the sense of [39] and conditional Rademacher averages [7]; see also [14] for other interesting applications.

COROLLARY 1. *Assume that $0 \leq Z - Z_i \leq 1$ for all $i = 1, \ldots, n$ and that for some constant $A \geq 1$,*

$$0 \leq \sum_{i=1}^{n} (Z - Z_i) \leq AZ.$$

*Then for all integers $q \geq 1$,*

$$(2.1) \qquad \|Z\|_q \leq \mathbb{E}[Z] + A\frac{q-1}{2},$$

*and for every real $q \geq 2$, then*

$$\|(Z - \mathbb{E}[Z])_+\|_q \leq \sqrt{\kappa}\left[\sqrt{Aq\mathbb{E}[Z]} + \frac{Aq}{2}\right]. \tag{2.2}$$

*Moreover, for all integers $q \geq 2$,*

$$\|(Z - \mathbb{E}[Z])_-\|_q \leq \sqrt{CqA\mathbb{E}[Z]},$$

*where $C < 1.131$.*

The next corollary provides a simple sub-Gaussian bound for the lower tail whenever $V^-$ is bounded by a nondecreasing function of $Z$. A similar phenomenon was observed in ([7], Theorem 6).

COROLLARY 2. *Assume that $V^- \leq g(Z)$ for some nondecreasing function $g$. Then for all integers $q \geq 2$,*

$$\|(Z - \mathbb{E}[Z])_-\|_q \leq \sqrt{Kq\mathbb{E}[g(Z)]}.$$

Finally, the following corollary of Theorem 3 deals with a generalization of self-bounded functionals that was already considered in [7].

COROLLARY 3. *Assume that $Z_i \leq Z$ for all $i = 1, \ldots, n$ and $V \leq WZ$ for a random variable $W \geq 0$. Then for all reals $q \geq 2$ and all $\theta \in (0, 1]$,*

$$\|Z\|_q \leq (1 + \theta)\mathbb{E}[Z] + \frac{\kappa}{2}\left(1 + \frac{1}{\theta}\right)q\|W\|_q.$$

*Also,*

$$\|(Z - \mathbb{E}[Z])_+\|_q \leq \sqrt{2\kappa q\|W\|_q\mathbb{E}[Z]} + \kappa q\|W\|_q.$$

*If $M$ denotes a positive random variable such that for every $1 \leq i \leq n$,*

$$0 \leq Z - Z_i \leq M,$$

*then we also have*

$$\|(Z - \mathbb{E}[Z])_-\|_q \leq \sqrt{C_2 q(\|M\|_q(2\mathbb{E}[Z] + 2q\|W\|_q) \vee q\|M\|_q^2)},$$

*where $C_2 < 2.42$ is as in Theorem 4.*

The proofs of Theorems 2–4 and of Corollaries 1–3 are developed in two steps. First, in Section 4, building on the modified $\phi$-Sobolev inequalities presented in Section 3, generalized Efron–Stein-type moment inequalities are established. These modified $\phi$-Sobolev/Efron–Stein inequalities play a role similar to the one played by modified log-Sobolev inequalities in the entropy method in [26, 27, 28] and [30]. Second, in Section 5 these general inequalities are used as main steps of an inductive proof of the main results. This second step may be regarded as an analogue of what is called in [28] the *Herbst argument* of the entropy method.

**3. Modified $\phi$-Sobolev inequalities.** The purpose of this section is to reveal some fundamental connections between $\phi$-entropies and modified $\phi$-Sobolev inequalities. The basic result is the duality formula of Lemma 1 implying the tensorization inequality which is at the basis of the modified $\phi$-Sobolev inequalities of Theorems 5 and 6. These theorems immediately imply the generalized Efron–Stein inequalities of Lemmas 3–5.

3.1. *$\phi$-entropies, duality and the tensorization property.* First we investigate so-called "tensorization" inequalities due to Latała and Oleszkiewicz [25] and Bobkov (see [26]). As of the time of writing this text, Chafaï [11] developed a framework for $\phi$-entropies and $\phi$-Sobolev inequalities.

We introduce some notation. Let $\mathbb{L}_1^+$ denote the convex set of nonnegative and integrable random variables $Z$. For any convex function $\phi$ on $\mathbb{R}_+$, let the $\phi$-entropy functional $H_\phi$ be defined for $Z \in \mathbb{L}_1^+$ by

$$H_\phi(Z) = \mathbb{E}[\phi(Z)] - \phi(\mathbb{E}[Z]).$$

Note that here and below we use the extended notion of expectation for a (not necessarily integrable) random variable $X$ defined as $\mathbb{E}[X] = \mathbb{E}[X_+] - \mathbb{E}[X_-]$ whenever either $X_+$ or $X_-$ is integrable.

The functional $H_\phi$ is said to satisfy the *tensorization property* if for every finite family $X_1, \ldots, X_n$ of independent random variables and every $(X_1, \ldots, X_n)$-measurable nonnegative and integrable random variable $Z$,

$$H_\phi(Z) \leq \sum_{i=1}^n \mathbb{E}[\mathbb{E}[\phi(Z)|X^{(i)}] - \phi(\mathbb{E}[Z|X^{(i)}])].$$

Observe that for $n = 2$ and setting $Z = g(X_1, X_2)$, the tensorization property reduces to the Jensen-type inequality

$$(3.1) \qquad H_\phi\left(\int g(x, X_2)\, d\mu_1(x)\right) \leq \int H_\phi(g(x, X_2))\, d\mu_1(x),$$

where $\mu_1$ denotes the distribution of $X_1$. Next we show that (3.1) implies the tensorization property. Indeed let $Y_1$ be distributed like $X_1$, and $Y_2$ be distributed like the $(n-1)$-tuple $X_2, \ldots, X_n$. Let $\mu_1$ and $\mu_2$ denote the corresponding distributions. The random variable $Z$ is a measurable function $g$ of the two independent random variables $Y_1$ and $Y_2$. By the Tonelli–Fubini theorem,

$$\begin{aligned}
H_\phi(Z) = \iint \Bigg( &\phi(g(y_1, y_2)) - \phi\left(\int g(y_1', y_2)\, d\mu_1(y_1')\right) \\
&+ \phi\left(\int g(y_1', y_2)\, d\mu_1(y_1)\right) \\
&- \phi\left(\iint g(y_1', y_2')\, d\mu_1(y_1')\, d\mu_2(y_2')\right) \Bigg) d\mu_1(y_1)\, d\mu_2(y_2)
\end{aligned}$$

$$= \int \left( \int \left[ \phi(g(y_1, y_2)) - \phi\left( \int g(y_1', y_2) \, d\mu_1(y_1') \right) \right] d\mu_1(y_1) \right) d\mu_2(y_2)$$

$$+ \int \left( \phi\left( \int g(y_1', y_2) \, d\mu_1(y_1') \right) \right.$$

$$\left. - \phi\left( \int \int g(y_1', y_2') \, d\mu_1(y_1') \, d\mu_2(y_2') \right) \right) d\mu_2(y_2)$$

$$= \int H_\phi(g(Y_1, y_2)) \, d\mu_2(y_2) + H_\phi\left( \int g(y_1', Y_2) \, d\mu_1(y_1') \right)$$

$$\leq \int H_\phi(g(Y_1, y_2)) \, d\mu_2(y_2) + \int H_\phi(g(y_1', Y_2)) \, d\mu_1(y_1'),$$

where the last step follows from the Jensen-type inequality (3.1).

If we turn back to the original notation, we get

$$H_\phi(Z) \leq \mathbb{E}[\mathbb{E}[\phi(Z)|X^{(1)}] - \phi(\mathbb{E}[Z|X^{(1)}])]$$

$$+ \int [H_\phi(Z(x_1, X_2, \ldots, X_n))] \, d\mu_1(x_1).$$

Proceeding by induction, (3.1) leads to the tensorization property for every $n$. We see that the tensorization property for $H_\phi$ is equivalent to what we could call the *Jensen property*, that is, (3.1) holds for every $\mu_1$, $X_2$ and $g$ such that $\int g(x, X_2) \, d\mu_1(x)$ is integrable.

Let $\Phi$ denote the class of functions $\phi$ which are continuous and convex on $\mathbb{R}_+$, twice differentiable on $\mathbb{R}_+^*$, and such that either $\phi$ is affine or $\phi''$ is strictly positive and $1/\phi''$ is concave.

It is shown in [25] (see also [26]) that there is a tight connection between the convexity of $H_\phi$ and the tensorization property. Also, $\phi \in \Phi$ implies the convexity of $H_\phi$; see [25]. However, this does not straightforwardly lead to Jensen's property when the distribution $\mu_1$ in (3.1) is not discrete. (See Appendix A.1 for an account of the consequences of the convexity of $\phi$-entropy.)

The easiest way to establish that for some function $\phi$ the functional $H_\phi$ satisfies the Jensen-like property is by following the lines of Ledoux's proof of the tensorization property for the "usual" entropy [which corresponds to the case $\phi(x) = x \log(x)$] and mimicking the duality argument used in one dimension to prove the usual Jensen's inequality, that is, to express $H_\phi$ as a supremum of affine functions.

Provided that $\phi \in \Phi$, our next purpose is to establish a duality formula for $\phi$-entropy of the form

$$H_\phi(Z) = \sup_{T \in \mathcal{T}} \mathbb{E}[\psi_1(T)Z + \psi_2(T)],$$

for convenient functions $\psi_1$ and $\psi_2$ on $\mathbb{R}_+$ and a suitable class of nonnegative variables $\mathcal{T}$. Such a formula obviously implies the convexity of $H_\phi$ but also Jensen's property and therefore the tensorization property for $H_\phi$. Indeed, considering again $Z$ as a function of $Y_1 = X_1$ and $Y_2 = (X_1, \ldots, Y_n)$ and assuming that a duality formula of the above form holds, we have

$$H_\phi\left( \int g(y_1, Y_2)\, d\mu_1(y_1) \right)$$

$$= \sup_{T \in \mathcal{T}} \int \left[ \psi_1(T(y_2)) \int g(y_1, y_2)\, d\mu_1(y_1) + \psi_2(T(y_2)) \right] d\mu_2(y_2)$$

$$\text{(by Fubini)}$$

$$= \sup_{T \in \mathcal{T}} \int \left( \int [\psi_1(T(y_2))g(y_1, y_2) + \psi_2(T(y_2))]\, d\mu_2(y_2) \right) d\mu_1(y_1)$$

$$\leq \int \left( \sup_{T \in \mathcal{T}} \int [\psi_1(T(y_2))g(y_1, y_2) + \psi_2(T(y_2))]\, d\mu_2(y_2) \right) d\mu_1(y_1)$$

$$= \int (H_\phi(g(y_1, Y_2)))\, d\mu_1(y_1).$$

LEMMA 1. *Let* $\phi \in \Phi$ *and* $Z \in \mathbb{L}_1^+$. *If* $\phi(Z)$ *is integrable, then*

$$H_\phi(Z) = \sup_{T \in \mathbb{L}_1^+, T \neq 0} \{ \mathbb{E}[(\phi'(T) - \phi'(\mathbb{E}[T]))(Z - T) + \phi(T)] - \phi(\mathbb{E}[T]) \}.$$

REMARK. This duality formula is almost identical to Proposition 4 in [11]. However, the proofs have different flavor. The proof given here is elementary.

PROOF. The case when $\phi$ is affine is trivial: $H_\phi$ equals zero, and so does the expression defined by the duality formula.

Note that the expression within the brackets on the right-hand side equals $H_\phi(Z)$ for $T = Z$, so the proof of Lemma 1 amounts to checking that

$$H_\phi(Z) \geq \mathbb{E}[(\phi'(T) - \phi'(\mathbb{E}[T]))(Z - T) + \phi(T)] - \phi(\mathbb{E}[T])$$

under the assumption that $\phi(Z)$ is integrable and $T \in \mathbb{L}_1^+$.

Assume first that $Z$ and $T$ are bounded and bounded away from 0. For any $\lambda \in [0, 1]$, we set $T_\lambda = (1 - \lambda)Z + \lambda T$ and

$$f(\lambda) = \mathbb{E}[(\phi'(T_\lambda) - \phi'(\mathbb{E}[T_\lambda]))(Z - T_\lambda)] + H_\phi(T_\lambda).$$

Our aim is to show that $f$ is nonincreasing on $[0, 1]$. Noticing that $Z - T_\lambda = \lambda(Z - T)$ and using our boundedness assumptions to differentiate under the expectation, we have

$$f'(\lambda) = -\lambda[\mathbb{E}[(Z - T)^2 \phi''(T_\lambda)] - (\mathbb{E}[Z - T])^2 \phi''(\mathbb{E}[T_\lambda])]$$

$$+ \mathbb{E}[(\phi'(T_\lambda) - \phi'(\mathbb{E}[T_\lambda]))(Z - T)]$$

$$+ \mathbb{E}[\phi'(T_\lambda)(T - Z)] - \phi'(\mathbb{E}[T_\lambda])\mathbb{E}[T - Z],$$

that is,

$$f'(\lambda) = -\lambda[\mathbb{E}[(Z-T)^2\phi''(T_\lambda)] - (\mathbb{E}[Z-T])^2\phi''(\mathbb{E}[T_\lambda])].$$

Now, by the Cauchy–Schwarz inequality,

$$(\mathbb{E}[Z-T])^2 = \left(\mathbb{E}\left[(Z-T)\sqrt{\phi''(T_\lambda)}\frac{1}{\sqrt{\phi''(T_\lambda)}}\right]\right)^2$$

$$\leq \mathbb{E}\left[\frac{1}{\phi''(T_\lambda)}\right]\mathbb{E}[(Z-T)^2\phi''(T_\lambda)].$$

Using the concavity of $1/\phi''$, Jensen's inequality implies that

$$\mathbb{E}\left[\frac{1}{\phi''(T_\lambda)}\right] \leq \frac{1}{\phi''(\mathbb{E}[T_\lambda])},$$

which leads to

$$(\mathbb{E}[Z-T])^2 \leq \frac{1}{\phi''(\mathbb{E}[T_\lambda])}\mathbb{E}[(Z-T)^2\phi''(T_\lambda)],$$

which is equivalent to $f'(\lambda) \leq 0$ and therefore $f(1) \leq f(0) = H_\phi(Z)$. This means that for any $T$, $\mathbb{E}[(\phi'(T) - \phi'(\mathbb{E}[T]))(Z-T)] + H_\phi(T) \leq H_\phi(Z)$.

In the general case we consider the sequences $Z_n = (Z \vee 1/n) \wedge n$ and $T_k = (T \vee 1/k) \wedge k$ and our purpose is to take the limit, as $k, n \to \infty$, in the inequality

$$H_\phi(Z_n) \geq \mathbb{E}[(\phi'(T_k) - \phi'(\mathbb{E}[T_k]))(Z_n - T_k) + \phi(T_k)] - \phi(\mathbb{E}[T_k]),$$

which we can also write as

$$(3.2) \quad \mathbb{E}[\psi(Z_n, T_k)] \geq -\phi'(\mathbb{E}[T_k])\mathbb{E}[Z_n - T_k] - \phi(\mathbb{E}[T_k]) + \phi(\mathbb{E}[Z_n]),$$

where $\psi(z,t) = \phi(z) - \phi(t) - (z-t)\phi'(t)$. Since we have to show that

$$(3.3) \quad \mathbb{E}[\psi(Z,T)] \geq -\phi'(\mathbb{E}[T])\mathbb{E}[Z-T] - \phi(\mathbb{E}[T]) + \phi(\mathbb{E}[Z])$$

with $\psi \geq 0$, we can always assume $[\psi(Z,T)]$ to be integrable [since otherwise (3.3) is trivially satisfied]. Taking the limit when $n$ and $k$ go to infinity on the right-hand side of (3.2) is easy, while the treatment of the left-hand side requires some care. Note that $\psi(z,t)$, as a function of $t$, decreases on $(0, z)$ and increases on $(z, +\infty)$. Similarly, as a function of $z$, $\psi(z,t)$ decreases on $(0, t)$ and increases on $(t, +\infty)$. Hence, for every $t$, $\psi(Z_n, t) \leq \psi(1, t) + \psi(Z, t)$, while for every $z$, $\psi(z, T_k) \leq \psi(z, 1) + \psi(z, T)$. Hence, given $k$,

$$\psi(Z_n, T_k) \leq \psi(1, T_k) + \psi(Z, T_k),$$

as $\psi((z \vee 1/n) \wedge n, T_k) \to \psi(z, T_k)$ for every $z$, we can apply the dominated convergence theorem to conclude that $\mathbb{E}[\psi(Z_n, T_k)]$ converges to $\mathbb{E}[\psi(Z, T_k)]$ as $n$ goes to infinity. Hence we have the following inequality:

$$(3.4) \quad \mathbb{E}[\psi(Z, T_k)] \geq -\phi'(\mathbb{E}[T_k])\mathbb{E}[Z - T_k] - \phi(\mathbb{E}[T_k]) + \phi(\mathbb{E}[Z]).$$

Now we also have $\psi(Z, T_k) \leq \psi(Z, 1) \psi(Z, T)$ and we can apply the dominated convergence theorem again to ensure that $\mathbb{E}[\psi(Z, T_k)]$ converges to $\mathbb{E}[\psi(Z, T)]$ as $k$ goes to infinity. Taking the limit as $k$ goes to infinity in (3.4) implies that (3.3) holds for every $T, Z \in \mathbb{L}_1^+$ such that $\phi(Z)$ is integrable and $\mathbb{E}[T] > 0$. If $Z \neq 0$ a.s., (3.3) is achieved for $T = Z$, while if $Z = 0$ a.s., it is achieved for $T = 1$ and the proof of the lemma is now complete in its full generality. $\quad \square$

REMARK. Note that since the supremum in the duality formula of Lemma 1 is achieved for $T = Z$ (or $T = 1$ if $Z = 0$), the duality formula remains true if the supremum is restricted to the class $\mathcal{T}_\phi$ of variables $T$ such that $\phi(T)$ is integrable. Hence the following alternative formula also holds:

$$(3.5) \qquad H_\phi(Z) = \sup_{T \in \mathcal{T}_\phi} \{ \mathbb{E}[(\phi'(T) - \phi'(\mathbb{E}[T]))(Z - T)] + H_\phi(T) \}.$$

REMARK. The duality formula of Lemma 1 takes the following (known) form for the "usual" entropy [which corresponds to $\phi(x) = x \log(x)$]:

$$\mathrm{Ent}(Z) = \sup_T \{ \mathbb{E}[(\log(T) - \log(\mathbb{E}[T]))Z] \},$$

where the supremum is extended to the set of nonnegative and integrable random variables $T$ with $\mathbb{E}[T] > 0$. Another case of interest is $\phi(x) = x^p$, where $p \in (1, 2]$. In this case, one has, by (3.5),

$$H_\phi(Z) = \sup_T \{ p \mathbb{E}[Z(T^{p-1} - (\mathbb{E}[T])^{p-1})] - (p-1) H_\phi(T) \},$$

where the supremum is extended to the set of nonnegative variables in $\mathbb{L}_p$.

REMARK. For the sake of simplicity we have focused on nonnegative variables and convex functions $\phi$ on $\mathbb{R}_+$. This restriction can be avoided and one may consider the case where $\phi$ is a convex function on $\mathbb{R}$ and define the $\phi$-entropy of a real-valued integrable random variable $Z$ by the same formula as in the nonnegative case. Assuming this time that $\phi$ is differentiable on $\mathbb{R}$ and twice differentiable on $\mathbb{R} \setminus \{0\}$, the proof of the duality formula above can be easily adapted to cover this case provided that $1/\phi''$ can be extended to a concave function on $\mathbb{R}$. In particular, if $\phi(x) = |x|^p$, where $p \in (1, 2]$, one gets

$$H_\phi(Z) = \sup_T \left\{ p \mathbb{E}\left[ Z\left( \frac{|T|^p}{T} - \frac{|\mathbb{E}[T]|^p}{\mathbb{E}[T]} \right) \right] - (p-1) H_\phi(T) \right\},$$

where the supremum is extended to $\mathbb{L}_p$. Note that for $p = 2$ this formula reduces to the classical one for the variance

$$\mathrm{Var}(Z) = \sup_T \{ 2 \mathrm{Cov}(Z, T) - \mathrm{Var}(T) \},$$

where the supremum is extended to the set of square integrable variables. This means that the tensorization inequality for the $\phi$-entropy also holds for convex functions $\phi$ on $\mathbb{R}$ under the condition that $1/\phi''$ is the restriction to $\mathbb{R} \setminus \{0\}$ of a concave function on $\mathbb{R}$.

3.2. *From $\phi$-entropies to $\phi$-Sobolev inequalities.*   Recall that our aim is to derive moment inequalities based on the tensorization property of $\phi$-entropy for an adequate choice of the function $\phi$ (namely, a properly chosen power function).

As a training example, we show how to derive the Efron–Stein inequality cited in Proposition 1 and a variant of it from the tensorization inequality of the variance, that is, the $\phi$-entropy when $\phi$ is defined on the whole real line as $\phi(x) = x^2$. Then

$$\mathrm{Var}(Z) \leq \mathbb{E}\left[\sum_{i=1}^{n} \mathbb{E}[(Z - \mathbb{E}[Z|X^{(i)}])^2 | X^{(i)}]\right]$$

and since conditionally on $X^{(i)}$, $Z_i'$ is an independent copy of $Z$, one has

$$\mathbb{E}[(Z - \mathbb{E}[Z|X^{(i)}])^2 | X^{(i)}] = \tfrac{1}{2}\mathbb{E}[(Z - Z_i')^2 | X^{(i)}],$$

which leads to Proposition 1. A useful variant may be obtained by noticing that $\mathbb{E}[Z|X^{(i)}]$ is the best $X^{(i)}$-measurable approximation of $Z$ in $\mathbb{L}_2$ which leads to

$$(3.6) \qquad \mathrm{Var}(Z) \leq \sum_{i=1}^{n} \mathbb{E}[(Z - Z_i)^2]$$

for any family of square integrable random variables $Z_i$'s such that $Z_i$ is $X^{(i)}$-measurable.

Next we generalize these symmetrization and variational arguments. The derivation of modified $\phi$-Sobolev inequalities will rely on the following properties of the elements of $\Phi$. The proofs of Proposition 2 and Lemma 2 are given in Appendix A.1.

PROPOSITION 2.   *If $\phi \in \Phi$, then both $\phi'$ and $x \to (\phi(x) - \phi(0))/x$ are concave functions on $(0, \infty)$.*

LEMMA 2.   *Let $\phi$ be a continuous and convex function on $\mathbb{R}_+$. Then, denoting by $\phi'$ the right derivative of $\phi$, for every $Z \in \mathbb{L}_1^+$, one has*

$$(3.7) \qquad H_\phi(Z) = \inf_{u \geq 0} \mathbb{E}[\phi(Z) - \phi(u) - (Z - u)\phi'(u)].$$

*Let $Z'$ be an independent copy of $Z$. Then*

$$(3.8) \qquad \begin{aligned} H_\phi(Z) &\leq \tfrac{1}{2}\mathbb{E}[(Z - Z')(\phi'(Z) - \phi'(Z'))] \\ &= \mathbb{E}[(Z - Z')_+(\phi'(Z) - \phi'(Z'))]. \end{aligned}$$

*If, moreover,* $\psi : x \to (\phi(x) - \phi(0))/x$ *is concave on* $\mathbb{R}_+^*$, *then*

$$(3.9) \qquad \begin{aligned} H_\phi(Z) &\leq \tfrac{1}{2}\mathbb{E}[(Z - Z')(\psi(Z) - \psi(Z'))] \\ &= \mathbb{E}[(Z - Z')_+(\psi(Z) - \psi(Z'))]. \end{aligned}$$

Note that by Proposition 2, we can apply (3.9) whenever $\phi \in \Phi$. In particular, for our target example where $\phi(x) = x^p$, with $p \in (1, 2]$, (3.9) improves on (3.8) within a factor $p$.

Modified $\phi$-Sobolev inequalities follow then from the tensorization inequality for $\phi$-entropy, the variational formula and the symmetrization inequality. The goal is to upper bound the $\phi$-entropy of a conveniently chosen convex function $f$ of the variable of interest $Z$. The results crucially depend on the monotonicity of the transformation $f$.

THEOREM 5. *Let* $X_1, \ldots, X_n$ *be independent random variables and let* $Z$ *be an* $(X_1, \ldots, X_n)$*-measurable random variable taking its values in an interval* $\mathcal{I}$. *Let* $V$, $V^+$ *and* $(Z_i)_{i \leq n}$ *be defined as in Section* 2.1.

*Let* $\phi \in \Phi$ *and let* $f$ *be a nondecreasing, nonnegative and differentiable convex function on* $\mathcal{I}$. *Let* $\psi$ *denote the function* $x \to (\phi(x) - \phi(0))/x$. *Then*

$$H_\phi(f(Z)) \leq \mathbb{E}[V^+ f'^2(Z)\psi'(f(Z))] \qquad \text{*if* } \psi \circ f \text{ *is convex.*}$$

*On the other hand, if* $(Z_i)_{i \leq n}$ *satisfy* $Z_i \leq Z$ *for all* $i \leq n$, *then*

$$H_\phi(f(Z)) \leq \tfrac{1}{2}\mathbb{E}[V f'^2(Z)\phi''(f(Z))] \qquad \text{*if* } \phi' \circ f \text{ *is convex.*}$$

PROOF. First fix $x < y$. Assume first that $g = \phi' \circ f$ is convex. We first check that

$$(3.10) \qquad \begin{aligned} \phi(f(y)) &- \phi(f(x)) - (f(y) - f(x))\phi'(f(x)) \\ &\leq \tfrac{1}{2}(y - x)^2 f'^2(y)\phi''(f(y)). \end{aligned}$$

Indeed, setting

$$h(t) = \phi(f(y)) - \phi(f(t)) - (f(y) - f(t))g(t),$$

we have

$$h'(t) = -g'(t)(f(y) - f(t)).$$

But for every $t \leq y$, the monotonicity and convexity assumptions on $f$ and $g$ yield

$$0 \leq g'(t) \leq g'(y) \quad \text{and} \quad 0 \leq f(y) - f(t) \leq (y - t)f'(y),$$

hence

$$-h'(t) \leq (y - t)f'(y)g'(y).$$

Integrating this inequality with respect to $t$ on $[x, y]$ leads to (3.10).

Under the assumption that $\psi \circ f$ is convex,

$$0 \leq f(y) - f(x) \leq (y - x) f'(y)$$

and

$$0 \leq \psi(f(y)) - \psi(f(x)) \leq (y - x) f'(y) \psi'(f(y)),$$

which leads to

$$(3.11) \quad (f(y) - f(x))(\psi(f(y)) - \psi(f(x))) \leq (x - y)^2 f'^2(y) \psi'(f(y)).$$

Now the tensorization inequality combined with the variational inequality (3.7) from Lemma 2 and (3.10) lead to

$$H_\phi(f(Z)) \leq \tfrac{1}{2} \sum_{i=1}^n \mathbb{E}[(Z - Z_i)^2 f'^2(Z) \phi''(f(Z))]$$

and therefore to the second inequality of the theorem.

The first inequality of the theorem follows in a similar way from inequality (3.9) and from (3.11).  $\square$

The case when $f$ is nonincreasing is handled by the following theorem.

THEOREM 6.  *Let $X_1, \ldots, X_n$ be independent random variables and let $Z$ be an $(X_1, \ldots, X_n)$-measurable random variable taking its values in some interval $\mathcal{I}$. Let $\phi \in \Phi$ and let $f$ be a nonnegative, nonincreasing and differentiable convex function on $\mathcal{I}$. Let $\psi$ denote the function $x \to (\phi(x) - \phi(0))/x$. For any random variable $\widetilde{Z} \leq \min_{1 \leq i \leq n} Z_i$,*

$$H_\phi(f(Z)) \leq \tfrac{1}{2} \mathbb{E}[V f'^2(\widetilde{Z}) \phi''(f(\widetilde{Z}))] \qquad \text{if } \phi' \circ f \text{ is convex,}$$

*while if $\psi \circ f$ is convex, we have*

$$H_\phi(f(Z)) \leq \mathbb{E}[V^+ f'^2(\widetilde{Z}) \psi'(f(\widetilde{Z}))]$$

*and*

$$H_\phi(f(Z)) \leq \mathbb{E}[V^- f'^2(Z) \psi'(f(Z))].$$

The proof of Theorem 6 parallels the proof of Theorem 5. It is included in Appendix A.1 for the sake of completeness.

REMARK.  As a first illustration, we may derive the modified logarithmic Sobolev inequalities in [7] using Theorems 5 and 6. Indeed, letting $f(z) = \exp(\lambda z)$ and $\phi(x) = x \log(x)$ leads to

$$H_\phi(f(Z)) \leq \lambda^2 \mathbb{E}[V^+ \exp(\lambda Z)],$$

if $\lambda \geq 0$, while if $\lambda \leq 0$, one has

$$H_\phi(f(Z)) \leq \lambda^2 \mathbb{E}[V^- \exp(\lambda Z)].$$

**4. Generalized Efron–Stein inequalities.** The purpose of this section is to prove the next three lemmas which relate different moments of $Z$ to $V$, $V^+$ and $V^-$. These lemmas are generalizations of the Efron–Stein inequality.

Recall the definitions of $(X_i), Z, (Z_i), (Z'_i), V^+, V^-, V$ and the constants $\kappa$ and $K$, given in Section 2.1.

LEMMA 3. *Let $q \geq 2$ be a real number and let $\alpha$ satisfy $q/2 \leq \alpha \leq q - 1$. Then*

$$\mathbb{E}[(Z - \mathbb{E}[Z])_+^q] \leq \mathbb{E}[(Z - \mathbb{E}[Z])_+^\alpha]^{q/\alpha} + \frac{q(q - \alpha)}{2} \mathbb{E}[V(Z - \mathbb{E}[Z])_+^{q-2}],$$

$$\mathbb{E}[(Z - \mathbb{E}[Z])_+^q] \leq \mathbb{E}[(Z - \mathbb{E}[Z])_+^\alpha]^{q/\alpha} + \alpha(q - \alpha)\mathbb{E}[V^+(Z - \mathbb{E}[Z])_+^{q-2}]$$

*and*

$$\mathbb{E}[(Z - \mathbb{E}[Z])_-^q] \leq \mathbb{E}[(Z - \mathbb{E}[Z])_-^\alpha]^{q/\alpha} + \alpha(q - \alpha)\mathbb{E}[V^-(Z - \mathbb{E}[Z])_-^{q-2}].$$

PROOF. Let $q$ and $\alpha$ be chosen in such a way that $1 \leq q/2 \leq \alpha \leq q - 1$. Let $\phi(x) = x^{q/\alpha}$. Applying Theorem 5 with $f(z) = (z - \mathbb{E}[Z])_+^\alpha$ leads to the first two inequalities. Finally, we may apply the third inequality of Theorem 6 with $f(z) = (z - \mathbb{E}[Z])_-^\alpha$ to obtain the third inequality of the lemma. □

The next lemma is a variant of Lemma 3 that may be convenient when dealing with positive random variables.

LEMMA 4. *Let $q$ denote a real number, $q \geq 2$ and $q/2 \leq \alpha \leq q - 1$. If for all $i = 1, \ldots, n$*

$$0 \leq Z_i \leq Z \qquad a.s.,$$

*then*

$$\mathbb{E}[Z^q] \leq \mathbb{E}[Z^\alpha]^{q/\alpha} + \frac{q(q - \alpha)}{2}\mathbb{E}[VZ^{q-2}].$$

PROOF. The lemma follows by choosing $q$ and $\alpha$ such that $1 \leq q/2 \leq \alpha \leq q - 1$, taking $\phi(x) = x^{q/\alpha}$ and applying Theorem 5 with $f(z) = z^\alpha$. □

The third lemma will prove useful when dealing with lower tails.

LEMMA 5. *If the increments $Z - Z_i$ or $Z - Z'_i$ are bounded by some positive random variable $M$, then*

(4.1)
$$\mathbb{E}[(Z - \mathbb{E}[Z])_-^q]$$
$$\leq \mathbb{E}[(Z - \mathbb{E}[Z])_-^\alpha]^{q/\alpha} + \frac{q(q - \alpha)}{2}\mathbb{E}[V(Z - \mathbb{E}[Z] - M)_-^{q-2}].$$

*If the increments $Z - Z'_i$ are bounded by some positive random variable $M$, then*

$$(4.2) \quad \begin{aligned} &\mathbb{E}[(Z - \mathbb{E}[Z])^q_-] \\ &\qquad \leq \mathbb{E}[(Z - \mathbb{E}[Z])^\alpha_-]^{q/\alpha} + \alpha(q - \alpha)\mathbb{E}[V^+(Z - \mathbb{E}[Z] - M)^{q-2}_-]. \end{aligned}$$

PROOF. If the increments $Z - Z_i$ or $Z - Z'_i$ are upper bounded by some positive random variable $M$, then we may also use the alternative bounds for the lower deviations stated in Theorem 6 to derive both inequalities. □

To obtain the main results of the paper, the inequalities of the lemmas above may be used by induction on the order of the moment. The details are worked out in the next section.

## 5. Proof of the main theorems.

We are now prepared to prove Theorems 1–3 and Corollaries 1 and 3.

To illustrate the method of proof on the simplest possible example, first we present the proof of Theorem 1. This proof relies on a technical lemma proved in Appendix A.2. Recall from Section 2.1 that $K$ is defined as $1/(e - \sqrt{e})$.

LEMMA 6. *For all integers $q \geq 4$, the sequence*

$$q \mapsto x_q = \left(\frac{q-1}{q}\right)^{q/2}\left(1 + \frac{1}{K}\left(\frac{q-2}{q-1}\right)^{(q-2)/2}\right)$$

*is bounded by 1. Also, $\lim_{q\to\infty} x_q = 1$.*

PROOF OF THEOREM 1. To prove the first inequality, assume that $V^+ \leq c$. Let $m_q$ be defined by

$$m_q = \|(Z - \mathbb{E}[Z])_+\|_q.$$

For $q \geq 3$, we obtain from the second inequality of Lemma 3, with $\alpha = q - 1$,

$$(5.1) \quad m^q_q \leq m^q_{q-1} + c(q-1)m^{q-2}_{q-2}.$$

Our aim is to prove that

$$(5.2) \quad m^q_q \leq (Kqc)^{q/2} \qquad \text{for } q \geq 2.$$

To this end, we proceed by induction. For $q = 2$, note that by the Efron–Stein inequality,

$$m^2_2 \leq \mathbb{E}[V^+] \leq c$$

and therefore (5.2) holds for $q = 2$.

Taking $q = 3$, since $m_1 \leq m_2 \leq \sqrt{c}$, we derive from (5.1) that

$$m_3^3 \leq 3c^{3/2}.$$

This implies that (5.2) also holds for $q = 3$.

Consider now $q \geq 4$ and assume that

$$m_j \leq \sqrt{Kjc}$$

for every $j \leq q - 1$. Then, it follows from (5.1) and two applications of the induction hypothesis that

$$m_q^q \leq K^{q/2} c^{q/2} \sqrt{q-1} (\sqrt{q-1})^{q-1} + \frac{K^{q/2}}{K} c^{q/2} (q-1) (\sqrt{q-2})^{q-2}$$

$$= (Kqc)^{q/2} \left( \left( \frac{q-1}{q} \right)^{q/2} + \frac{q-1}{Kq} \left( \frac{q-2}{q} \right)^{(q-2)/2} \right)$$

$$= (Kqc)^{q/2} \left( \frac{q-1}{q} \right)^{q/2} \left( 1 + \frac{1}{K} \left( \frac{q-2}{q-1} \right)^{(q-2)/2} \right).$$

The first part of the theorem then follows from Lemma 6.

To prove the second part, note that if, in addition, $V^- \leq c$, then applying the first inequality to $-Z$, we obtain

$$\| (Z - \mathbb{E}[Z])_- \|_q \leq K\sqrt{qc}.$$

The statement follows by noting that

$$\mathbb{E}[|Z - \mathbb{E}[Z]|^q] = \mathbb{E}[(Z - \mathbb{E}[Z])_+^q] + \mathbb{E}[(Z - \mathbb{E}[Z])_-^q] \leq 2(K\sqrt{qc})^q. \qquad \square$$

The proof of Theorems 2 and 3, given together below, is very similar to the proof of Theorem 1 above.

PROOF OF THEOREMS 2 AND 3.   It suffices to prove the first inequality of Theorems 2 and 3 since the second inequality of Theorem 2 follows from the first by replacing $Z$ by $-Z$.

We intend to prove by induction on $k$ that for all integers $k \geq 1$, all $q \in (k, k+1]$,

$$\| (Z - \mathbb{E}[Z])_+ \|_q \leq \sqrt{q \kappa_q c_q},$$

where either $c_q = \| V \|_{q/2 \vee 1}$ or $c_q = 2 \| V^+ \|_{q/2 \vee 1} (1 - 1/q)$.

For $k = 1$, it follows from Hölder's inequality, the Efron–Stein inequality and its variant (3.6) that

$$\| (Z - \mathbb{E}[Z])_+ \|_q \leq \sqrt{2 \| V^+ \|_1} \leq \sqrt{2 \kappa_q \| V^+ \|_{1 \vee q/2}}$$

and

$$\| (Z - \mathbb{E}[Z])_+ \|_q \leq \sqrt{\| V \|_{1 \vee q/2}} \leq \sqrt{\kappa_q \| V \|_{1 \vee q/2}}.$$

Assume the property holds for all integers smaller than some $k > 1$, and let us consider $q \in (k, k+1]$. Hölder's inequality implies that for every nonnegative random variable $Y$,

$$\mathbb{E}[Y(Z - \mathbb{E}[Z])_+^{q-2}] \leq \|Y\|_{q/2} \|(Z - \mathbb{E}[Z])_+\|_q^{q-2},$$

hence, using the first and second inequalities of Lemma 3 with $\alpha = q - 1$, we get

$$\|(Z - \mathbb{E}[Z])_+\|_q^q \leq \|(Z - \mathbb{E}[Z])_+\|_{q-1}^q + \frac{q}{2} c_q \|(Z - \mathbb{E}[Z])_+\|_q^{q-2}.$$

Defining

$$x_q = \|(Z - \mathbb{E}[Z])_+\|_q^q (q\kappa_q c_q)^{-q/2},$$

it suffices to prove that $x_q \leq 1$. With this notation the previous inequality becomes

$$x_q q^{q/2} c_q^{q/2} \kappa_q^{q/2} \leq x_{q-1}^{q/q-1} (q-1)^{q/2} c_{q-1}^{q/2} \kappa_{q-1}^{q/2} + \tfrac{1}{2} x_q^{1-2/q} q^{q/2} c_q^{q/2} \kappa_q^{q/2-1},$$

from which we derive, since $c_{q-1} \leq c_q$ and $\kappa_{q-1} \leq \kappa_q$,

$$x_q \leq x_{q-1}^{q/q-1} \left(1 - \frac{1}{q}\right)^{q/2} + \frac{1}{2\kappa_q} x_q^{1-2/q}.$$

Assuming, by induction, that $x_{q-1} \leq 1$, the previous inequality implies that

$$x_q \leq \left(1 - \frac{1}{q}\right)^{q/2} + \frac{1}{2\kappa_q} x_q^{1-2/q}.$$

Since the function

$$f_q : x \to \left(1 - \frac{1}{q}\right)^{q/2} + \frac{1}{2\kappa_q} x^{1-2/q} - x$$

is strictly concave on $\mathbb{R}_+$ and positive at $x = 0$, $f_q(1) = 0$ and $f_q(x_q) \geq 0$ imply that $x_q \leq 1$ as desired. $\square$

PROOF OF THEOREM 4. We use the notation $m_q = \|(Z - \mathbb{E}Z)_-\|_q$. For $a > 0$, the continuous function

$$x \to e^{-1/2} + \frac{1}{ax} e^{1/\sqrt{x}} - 1$$

decreases from $+\infty$ to $e^{-1/2} - 1 < 0$ on $(0, +\infty)$. Define $C_a$ as the unique zero of this function.

Since $C_1$ and $C_2$ are larger than $1/2$, it follows from Hölder's inequality, the Efron–Stein inequality and its variant (3.6) that for $q \in [1, 2]$,

$$\|(Z - \mathbb{E}[Z])_-\|_q \leq \sqrt{2\|V^+\|_1} \leq \sqrt{2\kappa_q \|V^+\|_{1 \vee q/2}}$$

and
$$\|(Z - \mathbb{E}[Z])_-\|_q \le \sqrt{\|V\|_{1 \vee q/2}} \le \sqrt{\kappa_q \|V\|_{1 \vee q/2}}.$$

In the rest of the proof the two cases may be dealt with together. The first case, belonging to the first assumption of Theorem 4, corresponds to $a = 1$, while the second corresponds to $a = 2$. Thus, we define

$$c_q = \begin{cases} \|V^+\|_{1 \vee q/2} \vee q \|M\|_q^2, & \text{when } a = 1, \\ \|V\|_{1 \vee q/2} \vee q \|M\|_q^2, & \text{when } a = 2. \end{cases}$$

For $q \ge 2$, either (4.2) or (4.1) with $\alpha = q - 1$ implies

$$(5.3) \qquad m_q^q \le m_{q-1}^q + q \mathbb{E}[V^+((Z - \mathbb{E}Z)_- + M)^{q-2}]$$

and

$$(5.4) \qquad m_q^q \le m_{q-1}^q + \frac{q}{2} \mathbb{E}[V((Z - \mathbb{E}Z)_- + M)^{q-2}].$$

We first deal with the case $q \in [2, 3)$. By the subadditivity of $x \to x^{q-2}$ for $q \in [2, 3]$, we have

$$((Z - \mathbb{E}Z)_- + M)^{q-2} \le M^{q-2} + (Z - \mathbb{E}[Z])_-^{q-2}.$$

Using Hölder's inequality, we obtain from (5.3) and (5.4) that

$$m_q^q \le m_{q-1}^q + q \|M\|_q^{q-2} \|V^+\|_{q/2} + q \|V^+\|_{q/2} m_q^{q-2}$$

and

$$m_q^q \le m_{q-1}^q + \frac{q}{2} \|M\|_q^{q-2} \|V\|_{q/2} + \frac{q}{2} \|V\|_{q/2} m_q^{q-2}.$$

Using the fact that $m_{q-1} \le \sqrt{c_{q-1}} \le \sqrt{c_q}$, those two latter inequalities imply

$$m_q^q \le c_q^{q/2} + \frac{q^{2-q/2}}{a} c_q^{q/2} + \frac{q}{a} c_q m_q^{q-2}.$$

Let $x_q = (\frac{m_q}{\sqrt{C_a q c_q}})^q$; then the preceding inequality translates into

$$x_q \le \left(\frac{1}{C_a c_q}\right)^{q/2} + \frac{1}{aC_a}((\sqrt{C_a} q)^{-q+2} + x_q^{1-2/q})$$

which in turn implies

$$x_q \le \frac{1}{2C_a} + \frac{1}{aC_a}(1 + x_q^{1-2/q})$$

since $q \ge 2$ and $C_a \ge 1$.

The function

$$g_q : x \to \frac{1}{2C_a} + \frac{1}{aC_a}(1 + x^{1-2/q}) - x$$

is strictly concave on $\mathbb{R}_+$ and positive at 0. Furthermore,

$$g_q(1) = \frac{4+a}{2aC_a} - 1 < 0,$$

since $C_a > (4+a)/2a$. Hence $g_q$ can be nonnegative at point $x_q$ only if $x_q \leq 1$, which settles the case $q \in [2,3]$.

We now turn to the case $q \geq 3$. We will prove by induction on $k \geq 2$ that for all $q \in [k, k+1)$, $m_q \leq \sqrt{qC_a\kappa_q c_q}$. By the convexity of $x \to x^{q-2}$ we have, for every $\theta \in (0,1)$,

$$((Z - \mathbb{E}Z)_- + M)^{q-2} = \left(\theta\frac{(Z - \mathbb{E}Z)_-}{\theta} + (1-\theta)\frac{M}{1-\theta}\right)^{q-2}$$
$$\leq \theta^{-q+3}M^{q-2} + (1-\theta)^{-q+3}(Z - \mathbb{E}[Z])_-^{q-2}.$$

Using Hölder's inequality, we obtain from (5.3) and (5.4) that

$$m_q^q \leq m_{q-1}^q + q\theta^{-q+3}\|M\|_q^{q-2}\|V^+\|_{q/2} + q(1-\theta)^{-q+3}\|V^+\|_{q/2}m_q^{q-2}$$

and

$$m_q^q \leq m_{q-1}^q + \frac{q}{2}\theta^{-q+3}\|M\|_q^{q-2}\|V\|_{q/2} + \frac{q}{2}(1-\theta)^{-q+3}\|V\|_{q/2}m_q^{q-2}.$$

Now assume by induction that $m_{q-1} \leq \sqrt{C_a(q-1)c_{q-1}}$. Since $c_{q-1} \leq c_q$, we have

$$m_q^q \leq C_a^{q/2}(q-1)^{q/2}c_q^{q/2} + \frac{1}{a}q^{-q+2}\theta^{-q+3}q^{q/2}c_q^{q/2} + \frac{1}{a}q(1-\theta)^{-q+3}c_q m_q^{q-2}.$$

Let $x_q = C_a^{-q/2}m_q^q(qc_q)^{-q/2}$. Then it suffices to show that $x_q \leq 1$ for all $q > 2$. Observe that

$$x_q \leq \left(1 - \frac{1}{q}\right)^{q/2} + \frac{1}{aC_a}(\theta^{-q+3}(\sqrt{C_a}q)^{-q+2} + (1-\theta)^{-q+3}x_q^{1-2/q}).$$

We choose $\theta$ minimizing

$$g(\theta) = \theta^{-q+3}(\sqrt{C_a}q)^{-q+2} + (1-\theta)^{-q+3},$$

that is, $\theta = 1/(\sqrt{C_a}q + 1)$. Since for this value of $\theta$

$$g(\theta) = \left(1 + \frac{1}{\sqrt{C_a}q}\right)^{q-2},$$

the bound on $x_q$ becomes

$$x_q \leq \left(1 - \frac{1}{q}\right)^{q/2} + \frac{1}{aC_a}\left(1 + \frac{1}{\sqrt{C_a}q}\right)^{q-2}\left(1 + \left(\frac{\sqrt{C_a}q}{1 + \sqrt{C_a}q}\right)(x_q^{1-2/q} - 1)\right).$$

Hence, using the elementary inequalities

$$\left(1 - \frac{1}{q}\right)^{q/2} \leq e^{-1/2} \quad \text{and} \quad \left(1 + \frac{1}{\sqrt{C_a}q}\right)^{q-2} \leq e^{1/\sqrt{C_a}},$$

we get

$$x_q \leq e^{-1/2} + \frac{e^{1/\sqrt{C_a}}}{aC_a}\left(\frac{\sqrt{C_a}q}{1 + \sqrt{C_a}q}\right)(x_q^{1-2/q} - 1).$$

Since the function

$$f_q: x \to e^{-1/2} + \frac{e^{1/\sqrt{C_a}}}{aC_a}\left(1 + \left(\frac{\sqrt{C_a}q}{1 + \sqrt{C_a}q}\right)(x^{1-2/q} - 1)\right) - x$$

is strictly concave on $\mathbb{R}_+$ and positive at 0 and $C_a$ is defined in such a way that $f_q(1) = 0$, $f_q$ can be nonnegative at $x_q$ only if $x_q \leq 1$, which proves the theorem by induction. $\quad\square$

PROOF OF COROLLARY 1. Applying Lemma 4 with $\alpha = q - 1$ leads to

$$\|Z\|_q^q \leq \|Z\|_{q-1}^q + \frac{q}{2}\mathbb{E}[VZ^{q-2}].$$

But by assumption, we have $V \leq AZ$, and therefore

$$\|Z\|_q^q \leq \|Z\|_{q-1}^q + \frac{qA}{2}\|Z\|_{q-1}^{q-1}$$

$$\leq \|Z\|_{q-1}^q\left[1 + \frac{qA}{2\|Z\|_{q-1}}\right].$$

Since for any nonnegative real number $u$, $1 + uq \leq (1+u)^q$ for $u \geq 0$,

$$\|Z\|_q^q \leq \|Z\|_{q-1}^q\left(1 + \frac{A}{2\|Z\|_{q-1}}\right)^q$$

or, equivalently,

$$\|Z\|_q \leq \|Z\|_{q-1} + \frac{A}{2}.$$

Thus, $\|Z\|_q \leq \|Z\|_1 + (A/2)(q-1)$ by induction, and (2.1) follows.

To prove (2.2), note first that by Theorem 3,

$$\|(Z - \mathbb{E}Z)_+\|_q \leq \sqrt{\kappa q\|V\|_{q/2}} \leq \sqrt{\kappa qA\|Z\|_{q/2}}.$$

Let $s$ be the smallest integer such that $q/2 \leq s$. Then (2.1) yields

$$\|Z\|_{q/2} \leq \mathbb{E}[Z] + \frac{A(s-1)}{2} \leq \mathbb{E}[Z] + \frac{Aq}{4}$$

so that

$$\|(Z - \mathbb{E}[Z])_+\|_q \leq \sqrt{\kappa}\left[\sqrt{qA\mathbb{E}[Z] + \frac{q^2A^2}{4}}\right]$$

$$\leq \sqrt{\kappa}\left[\sqrt{qA\mathbb{E}[Z]} + \frac{qA}{2}\right]$$

and inequality (2.2) follows.

In order to prove the last inequality of Corollary 1, we first define $C$ as the unique positive root of the equation

$$e^{-1/2} + \frac{1}{2C}e^{-1+1/C} - 1 = 0.$$

We derive from the upper bound $V \leq AZ$ and the modified Efron–Stein inequality (3.6) that

$$(\mathbb{E}|Z - \mathbb{E}Z|)^2 \leq \mathbb{E}[(Z - \mathbb{E}Z)^2] \leq A\mathbb{E}Z.$$

Since $C > 1$, this proves the inequality for $q = 1$ and $q = 2$. For $q \geq 3$, we assume, by induction, that $m_k \leq \sqrt{CkA\mathbb{E}[Z]}$ for $k = q-2$ and $k = q-1$ and use $V \leq AZ$ together with (4.1) with $\alpha = q - 1$. This gives

$$m_q^q \leq m_{q-1}^q + \frac{q}{2}A\mathbb{E}[Z((Z - \mathbb{E}[Z])_- + 1)^{q-2}].$$

Recall Chebyshev's negative association inequality which asserts that if $f$ is nondecreasing and $g$ is nonincreasing, then

$$\mathbb{E}[fg] \leq \mathbb{E}[f]\mathbb{E}[g].$$

Since the function $z \to ((z - \mathbb{E}[Z])_- + 1)^{q-2}$ decreases, by Chebyshev's negative association inequality, the previous inequality implies

$$m_q^q \leq m_{q-1}^q + \frac{q}{2}A\mathbb{E}[Z]\mathbb{E}[((Z - \mathbb{E}Z)_- + 1)^{q-2}].$$

Thus, this inequality becomes

$$m_q^q \leq m_{q-1}^q + \frac{q}{2}A\mathbb{E}[Z](1 + m_{q-2})^{q-2}$$

and therefore our induction assumption yields

$$m_q^q \leq \left(1 - \frac{1}{q}\right)^{q/2}(CqA\mathbb{E}[Z])^{q/2} + \frac{(CqA\mathbb{E}[Z])^{q/2}}{2C}\left[\frac{1}{\sqrt{CqA\mathbb{E}[Z]}} + \sqrt{1 - \frac{2}{q}}\right]^{q-2}.$$

Now we use the fact that since $Z$ is nonnegative, $m_q \leq \mathbb{E}Z$. Then we may always assume that $CqA \leq \mathbb{E}Z$, since otherwise the last inequality of Corollary 1 is implied by this crude upper bound. Combining this inequality with $A > 1$ leads to

$$\frac{1}{\sqrt{CqA\mathbb{E}[Z]}} \leq \frac{1}{Cq},$$

so that plugging this in the inequality above and setting $x_q = m_q^q (CqA\mathbb{E}Z)^{-q/2}$, we derive that

$$x_q \le \left(1 - \frac{1}{q}\right)^{q/2} + \frac{1}{2C}\left(\frac{1}{Cq} + \sqrt{1 - \frac{2}{q}}\right)^{q-2}.$$

Now we claim that

$$(5.5) \qquad \left(\frac{1}{Cq} + \sqrt{1 - \frac{2}{q}}\right)^{q-2} \le e^{-1+1/C}.$$

Indeed, (5.5) may be checked numerically for $q = 3$, while for $q \ge 4$, combining

$$\sqrt{1 - \frac{2}{q}} \le 1 - \frac{1}{q} - \frac{1}{2q^2}$$

with $\ln(1 + u) \le u$ leads to

$$\ln\left[\left(\frac{1}{Cq} + \sqrt{1 - \frac{2}{q}}\right)^{q-2}\right] \le -1 + \frac{1}{C} + \frac{1}{q}\left(\frac{3}{2} + \frac{1}{q} - \frac{2}{C}\right)$$

$$\le -1 + \frac{1}{C} + \frac{1}{q}\left(\frac{7}{4} - \frac{2}{C}\right),$$

which, since $C < 8/7$, implies (5.5). Hence

$$x_q \le \left(1 - \frac{1}{q}\right)^{q/2} + \frac{1}{2C}e^{-1+1/C} \le e^{-1/2} + \frac{1}{2C}e^{-1+1/C},$$

which, by definition of $C$, means that $x_q \le 1$, completing the proof of the third inequality.  $\square$

PROOF OF COROLLARY 2.  This corollary follows by noting that if $V^-$ is bounded by a nondecreasing function of $Z$, then by negative association,

$$\mathbb{E}[V^-(Z - \mathbb{E}[Z])_-^{q-2}]$$

$$\le \mathbb{E}[g(Z)(Z - \mathbb{E}[Z])_-^{q-2}]$$

$$\le \mathbb{E}[g(Z)]\mathbb{E}[(Z - \mathbb{E}[Z])_-^{q-2}].$$

Thus, writing

$$m_q = \|(Z - \mathbb{E}[Z])_-\|_q,$$

we have

$$m_q^q \le m_{q-1}^q + \mathbb{E}[g(Z)](q-1)m_{q-2}^{q-2}.$$

This recursion is identical to the one appearing in the proof of Theorem 1, so the rest of the proof is identical to that of Theorem 1. □

PROOF OF COROLLARY 3.  Let $q$ be a number, $q \geq 2$. Let $\theta > 0$. Then

$$\|(Z - \mathbb{E}[Z])_+\|_q \leq \sqrt{\kappa q \|WZ\|_{q/2}} \qquad \text{(by Theorem 3)}$$

$$\leq \sqrt{\kappa q \|Z\|_q \|W\|_q} \qquad \text{(by Hölder's inequality)}$$

$$\leq \frac{\theta}{2}\|Z\|_q + \frac{\kappa q}{2\theta}\|W\|_q$$

$$[\text{for } \theta > 0, \text{since } \sqrt{ab} \leq (a^2 + b^2)/2 \text{ for } a, b \geq 0].$$

Now $Z \geq 0$ implies that $\|(Z - \mathbb{E}[Z])_-\|_q \leq \mathbb{E}[Z]$ and we have $\|Z\|_q \leq \mathbb{E}[Z] + \|(Z - \mathbb{E}[Z])_+\|_q$. Hence, for $0 < \theta \leq 1$,

$$\|Z\|_q \leq \frac{1}{1 - \theta/2}\mathbb{E}[Z] + \frac{\kappa q}{2\theta(1 - \theta/2)}\|W\|_q$$

$$\leq (1 + \theta)\mathbb{E}[Z] + \frac{\kappa q}{2}\left(1 + \frac{1}{\theta}\right)\|W\|_q,$$

concluding the proof of the first statement. To prove the second inequality, note that

$$\|(Z - \mathbb{E}[Z])_+\|_q$$

$$\leq \sqrt{\kappa q \|WZ\|_{q/2}} \qquad \text{(by Theorem 3)}$$

$$\leq \sqrt{\kappa q \|W\|_q \|Z\|_q} \qquad \text{(by Hölder's inequality)}$$

$$\leq \sqrt{\kappa q \|W\|_q (2\mathbb{E}[Z] + \kappa q \|W\|_q)} \qquad \text{(by the first inequality with } \theta = 1)$$

$$\leq \sqrt{2\kappa q \|W\|_q \mathbb{E}[Z]} + \kappa q \|W\|_q,$$

as desired. □

**6. Sums of random variables.**  In this section we show how the results stated in Section 2 imply some classical moment inequalities for sums of independent random variables such as the Khinchine–Kahane, Marcinkiewicz and Rosenthal inequalities. In all cases, the proof basically does not require any further work. Also, we obtain explicit constants which only depend on $q$. These constants are not optimal, though in some cases their dependence on $q$ is of the right order. For more information on these and related inequalities we refer to [13].

The simplest example is the case of Khinchine's inequality:

THEOREM 7 (Khinchine's inequality). *Let $a_1, \ldots, a_n$ be nonnegative constants, and let $X_1, \ldots, X_n$ be independent Rademacher variables (i.e., with $\mathbb{P}\{X_i = -1\} = \mathbb{P}\{X_i = 1\} = 1/2$). If $Z = \sum_{i=1}^{n} a_i X_i$, then for any integer $q \geq 2$,*

$$\|(Z)_+\|_q = \|(Z)_-\|_q \leq \sqrt{2Kq}\sqrt{\sum_{i=1}^{n} a_i^2}$$

*and*

$$\|Z\|_q \leq 2^{1/q}\sqrt{2Kq}\sqrt{\sum_{i=1}^{n} a_i^2},$$

*where $K = 1/(e - \sqrt{e}) < 0.935$.*

PROOF. We may use Theorem 1. Since

$$V^+ = \sum_{i=1}^{n} \mathbb{E}[(a_i(X_i - X_i'))_+^2 \,|\, X_i] = 2\sum_{i=1}^{n} a_i^2 \mathbb{1}_{a_i X_i > 0} \leq 2\sum_{i=1}^{n} a_i^2,$$

the result follows. □

Note also that using a symmetrization argument (see, e.g., [13], Lemma 1.2.6), Khinchine's inequality above implies Marcinkiewicz's inequality: if $X_1, \ldots, X_n$ are independent centered random variables, then for any $q \geq 2$,

$$\left\|\sum_{i=1}^{n} X_i\right\|_q \leq 2^{1+1/q}\sqrt{2Kq}\sqrt{\left\|\sum_i X_i^2\right\|_{q/2}}.$$

The next two results are Rosenthal-type inequalities for sums of independent nonnegative and centered random variables. The following inequality is very similar to inequality $(H_r)$ in [16] which follows from an improved Hoffmann–Jørgensen inequality of [24]. Note again that we obtain the result without further work.

THEOREM 8. *Define*

$$Z = \sum_{i=1}^{n} X_i,$$

*where $X_i$ are independent and nonnegative random variables. Then for all integers $q \geq 1$ and $\theta \in (0, 1)$,*

$$\|(Z - \mathbb{E}[Z])_+\|_q \leq \sqrt{2\kappa q \left\|\max_i |X_i|\right\|_q \mathbb{E}[Z]} + \kappa q \left\|\max_i |X_i|\right\|_q,$$

$$\|(Z - \mathbb{E}[Z])_-\|_q \leq \sqrt{Kq\sum_i \mathbb{E}[X_i^2]}$$

*and*

$$\|Z\|_q \leq (1 + \theta)\mathbb{E}[Z] + \frac{\kappa}{2}q\left(1 + \frac{1}{\theta}\right)\left\|\max_{1 \leq i \leq n} X_i\right\|_q.$$

PROOF.   We may use Corollary 3 to get the first and the third inequalities; just note that

$$V = \sum_{i=1}^{n} X_i^2 \leq WZ,$$

where

$$W = \max_{1 \leq i \leq n} X_i.$$

In order to get the second inequality, just observe that

$$V^- \leq \sum_i \mathbb{E}[{X'_i}^2],$$

and apply Theorem 1 to $-Z$.   □

Next we use the previous result to derive a Rosenthal-type inequality for sums of centered variables. In spite of the simplicity of the proof, the dependence of the constants on $q$ matches the best known bounds. (See [35] which extends the theorem below for martingales.)

THEOREM 9.   *Let* $X_i$, $i = 1, \ldots, n$, *be independent centered random variables. Define*

$$Z = \sum_{i=1}^{n} X_i, \qquad \sigma^2 = \sum_i \mathbb{E}[X_i^2], \qquad Y = \max_{1 \leq i \leq n} |X_i|.$$

*Then for any integer* $q \geq 2$ *and* $\theta \in (0, 1)$,

$$\|(Z)_+\|_q \leq \sigma\sqrt{2\kappa(2 + \theta)q} + q\kappa\sqrt{1 + \frac{1}{\theta}}\|Y\|_q.$$

PROOF.   We use Theorem 2. Note that

$$V^+ = \sum_i X_i^2 + \sum_i \mathbb{E}[{X'_i}^2].$$

Thus,

$$\|(Z)_+\|_q \leq \sqrt{2\kappa q\|V^+\|_{q/2}} \qquad \text{(by Theorem 2)},$$

$$\leq \sqrt{2\kappa q}\sqrt{\left(\sum_i \mathbb{E}[{X'_i}^2]\right) + \left\|\sum_i X_i^2\right\|_{q/2}}$$

$$\leq \sqrt{2\kappa q} \sqrt{\sum_i \mathbb{E}[{X_i'}^2] + (1+\theta)\sum_i \mathbb{E}[X_i^2] + \frac{\kappa q}{2}\left(1+\frac{1}{\theta}\right)\|Y^2\|_{q/2}}$$

$$\text{(by Theorem 8)}$$

$$= \sqrt{2\kappa q} \sqrt{(2+\theta)\sum_i \mathbb{E}[X_i^2] + \frac{\kappa q}{2}\left(1+\frac{1}{\theta}\right)\|Y^2\|_{q/2}}$$

$$\leq \sigma\sqrt{2\kappa(2+\theta)q} + q\kappa\sqrt{1+\frac{1}{\theta}}\|Y\|_q. \qquad \square$$

**7. Suprema of empirical processes.** In this section we apply the results of Section 2 to derive moment bounds for suprema of empirical processes. In particular, the main result of this section, Theorem 12, may be regarded as an analogue of Talagrand's inequality [40] for moments. Indeed, Talagrand's exponential inequality may be easily deduced from Theorem 12 by bounding the moment generating function by bounding all moments.

As a first illustration, we point out that the proof of Khinchine's inequality in the previous section extends, in a straightforward way, to an analogous supremum:

THEOREM 10. *Let $\mathcal{T} \subset \mathbb{R}^n$ be a set of vectors $t = (t_1, \ldots, t_n)$ and let $X_1, \ldots, X_n$ be independent Rademacher variables. If $Z = \sup_{t \in \mathcal{T}} \sum_{i=1}^n t_i X_i$, then for any integer $q \geq 2$,*

$$\|(Z - \mathbb{E}[Z])_+\|_q \leq \sqrt{2Kq} \sup_{t \in \mathcal{T}} \sqrt{\sum_{i=1}^n t_i^2},$$

*where $K = 1/(e - \sqrt{e}) < 0.935$, and*

$$\|(Z - \mathbb{E}[Z])_-\|_q \leq \sqrt{2C_1 q} \sup_{t \in \mathcal{T}} \sqrt{\sum_{i=1}^n t_i^2} \vee 2\sqrt{C_1}q \sup_{i,t}|t_i|,$$

*where $C_1$ is defined as in Theorem 4.*

Before stating the main result of the section, we mention the following consequence of Corollary 3.

THEOREM 11. *Let $\mathcal{F}$ be a countable class of nonnegative functions defined on some measurable set $\mathcal{X}$. Let $X_1, \ldots, X_n$ denote a collection of $\mathcal{X}$-valued independent random variables. Let $Z = \sup_{f \in \mathcal{F}} \sum_i f(X_i)$ and*

$$M = \max_{1 \leq i \leq n} \sup_{f \in \mathcal{F}} f(X_i).$$

*Then, for all $q \geq 2$ and $\theta \in (0, 2)$,*

$$\|Z\|_q \leq (1 + \theta)\mathbb{E}[Z] + \frac{\kappa}{2}q\left(1 + \frac{1}{\theta}\right)\|M\|_q.$$

Next we introduce the relevant quantities for the statement and proof of our main theorem about moments of centered empirical processes.

Let $\mathcal{F}$ denote a countable class of measurable functions from $\mathcal{X} \to \mathbb{R}$. Let $X_1, \ldots, X_n$ denote independent $\mathcal{X}$-valued random variables such that for all $f \in \mathcal{F}$ and $i = 1, \ldots, n$, $\mathbb{E}f(X_i) = 0$. Let

$$Z = \sup_{f \in \mathcal{F}}\left|\sum_{i=1}^{n} f(X_i)\right|.$$

The fluctuations of an empirical process are known to be characterized by two quantities that coincide when the process is indexed by a singleton. The strong variance $\Sigma^2$ is defined as

$$\Sigma^2 = \mathbb{E}\left[\sup_f \sum_i f^2(X_i)\right],$$

while the weak variance $\sigma^2$ is defined by

$$\sigma^2 = \sup_f \mathbb{E}\left[\sum_i f^2(X_i)\right].$$

A third quantity appearing in the moment and and tail bounds is

$$M = \sup_{i,f} |f(X_i)|.$$

Before stating the main theorem, we first establish a connection between the weak and the strong variances of an empirical process:

LEMMA 7.

$$\Sigma^2 \leq \sigma^2 + 32\sqrt{\mathbb{E}[M^2]}\mathbb{E}[Z] + 8\mathbb{E}[M^2].$$

If the functions in $\mathcal{F}$ are uniformly bounded, then $\Sigma$ may be upper bounded by a quantity that depends on $\sigma$ and $\mathbb{E}[Z]$ thanks to the contraction principle (see [30]). Giné, Latała and Zinn [16] combine the contraction principle with a Hoffmann–Jørgensen-type inequality. To follow their reasoning, we need the following lemma.

LEMMA 8. *Let $\varepsilon_1, \ldots, \varepsilon_n$ denote independent Rademacher variables. Let $\lambda > 4$ and define $t_0 = \sqrt{\lambda \mathbb{E}[M^2]}$. Then*

$$\mathbb{E}\left[\sup_f \left|\sum_i \varepsilon_i f^2(X_i) \mathbb{1}_{\sup_f |f(X_i)| > t_0}\right|\right] \leq \frac{1}{(1 - 2/\sqrt{\lambda})^2}\mathbb{E}[M^2].$$

The proof of this lemma is postponed to Appendix A.3.

PROOF OF LEMMA 7. Let $\varepsilon_1, \ldots, \varepsilon_n$ denote independent Rademacher random variables, and let $t_0 = \sqrt{\lambda \mathbb{E}[M^2]}$. Then

$$\Sigma^2 \leq \mathbb{E}\left[\sup_f \left| \sum_i f^2(X_i) - \mathbb{E}[f^2(X_i)] \right| \right] + \sup_f \mathbb{E}\left[ \sum_i f^2(X_i) \right]$$

$$\leq \sigma^2 + 2\mathbb{E}\left[ \sup_f \left| \sum_i \varepsilon_i f^2(X_i) \right| \right]$$

(by the symmetrization inequalities [29], Lemma 6.3)

$$\leq \sigma^2 + 2\mathbb{E}\left[ \sup_f \left| \sum_i \varepsilon_i f^2(X_i) \mathbb{1}_{\sup_f |f(X_i)| \leq t_0} \right| \right]$$

$$+ 2\mathbb{E}\left[ \sup_f \left| \sum_i \varepsilon_i f^2(X_i) \mathbb{1}_{\sup_f |f(X_i)| > t_0} \right| \right]$$

$$\leq \sigma^2 + 4t_0 \mathbb{E}\left[ \sup_f \left| \sum_i \varepsilon_i f(X_i) \right| \right] + 2\mathbb{E}\left[ \sup_f \left| \sum_i \varepsilon_i f^2(X_i) \mathbb{1}_{\sup_f |f(X_i)| > t_0} \right| \right]$$

(the contraction principle for Rademacher averages [29], Lemma 6.5

since $u \mapsto u^2/(2t_0)$ is contracting on $[-t_0, t_0]$)

$$\leq \sigma^2 + 4t_0 \mathbb{E}\left[ \sup_f \left| \sum_i \varepsilon_i f(X_i) \right| \right] + \frac{2}{(1 - 2/\sqrt{\lambda})^2} \mathbb{E}[M^2] \qquad \text{(by Lemma 8)}$$

$$\leq \sigma^2 + 8\sqrt{\lambda \mathbb{E}[M^2]} \|Z\|_1 + \frac{2}{(1 - 2/\sqrt{\lambda})^2} \mathbb{E}[M^2]$$

which, by taking $\lambda = 16$, completes the proof. $\square$

The next theorem offers two upper bounds for the moments of suprema of centered empirical processes. The first inequality improves inequality (3) of [35]. The second inequality is a version of Proposition 3.1 of [16]. It follows from the first combined with Lemma 8.

THEOREM 12. *Let $\mathcal{F}$ denote a countable class of measurable functions from $\mathcal{X} \to \mathbb{R}$. Let $X_1, \ldots, X_n$ denote independent $\mathcal{X}$-valued random variables such that for all $f \in \mathcal{F}$ and $i = 1, \ldots, n$, $\mathbb{E}f(X_i) = 0$. Let*

$$Z = \sup_{f \in \mathcal{F}} \left| \sum_{i=1}^n f(X_i) \right|.$$

*Then for all $q \geq 2$,*

$$\|(Z - \mathbb{E}[Z])_+\|_q \leq \sqrt{2\kappa q}(\Sigma + \sigma) + 2\kappa q\left(\|M\|_q + \sup_{i,f \in \mathcal{F}} \|f(X_i)\|_2\right),$$

*and furthermore,*

$$\|Z\|_q \leq 2\mathbb{E}Z + 2\sigma\sqrt{2\kappa q} + 20\kappa q\|M\|_q + 4\sqrt{\kappa q}\|M\|_2.$$

PROOF.  The proof uses Theorem 2 which states that

$$\|(Z - \mathbb{E}[Z])_+\|_q \leq \sqrt{2\kappa q\|V^+\|_{q/2}}.$$

We may bound $V^+$ as follows:

$$V^+ \leq \sup_{f \in \mathcal{F}} \sum_{i=1}^n \mathbb{E}[(f(X_i) - f(X_i'))^2 | X_1^n]$$

$$\leq \sup_{f \in \mathcal{F}} \sum_{i=1}^n (\mathbb{E}[f(X_i)^2] + f(X_i)^2)$$

$$\leq \sup_{f \in \mathcal{F}} \sum_{i=1}^n \mathbb{E}[f(X_i)^2] + \sup_{f \in \mathcal{F}} \sum_{i=1}^n f(X_i)^2.$$

Thus, by Minkowski's inequality and the Cauchy–Schwarz inequality,

$$\sqrt{\|V^+\|_{q/2}} \leq \sqrt{\sup_{f \in \mathcal{F}} \sum_{i=1}^n \mathbb{E}[f(X_i)^2] + \left\|\sup_{f \in \mathcal{F}} \sum_{i=1}^n f(X_i)^2\right\|_{q/2}}$$

$$\leq \sigma + \left\|\sup_{f \in \mathcal{F}} \sqrt{\sum_{i=1}^n f(X_i)^2}\right\|_q$$

$$= \sigma + \left\|\sup_{f \in \mathcal{F}} \sup_{\alpha \,:\, \|\alpha\|_2 \leq 1} \sum_{i=1}^n \alpha_i f(X_i)\right\|_q$$

$$\leq \sigma + \Sigma + \left\|\left(\sup_{f \in \mathcal{F}, \alpha \,:\, \|\alpha\|_2 \leq 1} \sum_{i=1}^n \alpha_i f(X_i)\right.\right.$$
$$\left.\left. - \mathbb{E}\left[\sup_{f \in \mathcal{F}, \alpha \,:\, \|\alpha\|_2 \leq 1} \sum_{i=1}^n \alpha_i f(X_i)\right]\right)_+\right\|_q.$$

The last summand may be upper bounded again by Theorem 2. Indeed, the corresponding $V^+$ is not more than

$$\max_i \sup_{f \in \mathcal{F}} f^2(X_i) + \max_i \sup_{f \in \mathcal{F}} \mathbb{E}[f^2(X_i)],$$

and thus

$$\left\| \left( \sup_{f \in \mathcal{F}, \alpha : \|\alpha\|_2 \le 1} \sum_{i=1}^n \alpha_i f(X_i) - \mathbb{E}\left[ \sup_{f \in \mathcal{F}, \alpha : \|\alpha\|_2 \le 1} \sum_{i=1}^n \alpha_i f(X_i) \right] \right)_+ \right\|_q$$

$$\le \sqrt{2\kappa q}\left( \|M\|_q + \max_i \sup_{f \in \mathcal{F}} \|f(X_i)\|_2 \right).$$

This completes the proof of the first inequality of the theorem. The second inequality follows because by nonnegativity of $Z$, $\|(Z - \mathbb{E}[Z])_-\|_q \le \mathbb{E}Z$ and therefore $\|Z\|_q \le \mathbb{E}Z + \|(Z - \mathbb{E}[Z])_+\|_q$ and since by the first inequality, combined with Lemma 7, we have

$$\|(Z - \mathbb{E}[Z])_+\|_q \le \sqrt{2\kappa q}\left( \sigma + \sqrt{32\sqrt{\mathbb{E}[M^2]\mathbb{E}[Z]}} + \sqrt{8\mathbb{E}[M^2]} + \sigma \right)$$

$$+ 2\kappa q \left( \|M\|_q + \sup_{i,f \in \mathcal{F}} \|f(X_i)\|_2 \right)$$

$$\le \mathbb{E}[Z] + 2\sigma\sqrt{2\kappa q} + 16\kappa\sqrt{\mathbb{E}[M^2]} + \sqrt{16\kappa q \mathbb{E}[M^2]}$$

$$+ 2\kappa q \left( \|M\|_q + \sup_{i,f \in \mathcal{F}} \|f(X_i)\|_2 \right)$$

(using the inequality $\sqrt{ab} \le a + b/4$).

Using $\|M\|_2 \le \|M\|_q$ and $\sup_{i,f \in \mathcal{F}} \|f(X_i)\|_2 \le \|M\|_2$, we obtain the desired result. $\square$

## 8. Conditional Rademacher averages.
Let $\mathcal{F}$ be a countable class of measurable real-valued functions. The conditional Rademacher average is defined by

$$Z = \mathbb{E}\left[ \sup_{f \in \mathcal{F}} \left| \sum_i \varepsilon_i f(X_i) \right| \,\middle|\, X_1^n \right],$$

where the $\varepsilon_i$ are i.i.d. Rademacher random variables. Conditional Rademacher averages play a distinguished role in probability in Banach spaces and in statistical learning theory (see, e.g., [1, 2, 3, 22, 23]). When the set of functions is bounded, $Z$ has been shown to satisfy a Bernstein-like inequality [7]. Here we provide bounds on the growth of moments in the general case.

THEOREM 13. *Let $Z$ denote a conditional Rademacher average and let* $M = \sup_{i,f} f(X_i)$. *Then*

$$\|(Z - \mathbb{E}[Z])_+\|_q \le \sqrt{2\kappa q \|M\|_q \mathbb{E}[Z]} + \kappa q \|M\|_q$$

*and*

$$\|(Z - \mathbb{E}[Z])_-\|_q \le \sqrt{2C_2}\{\sqrt{q \|M\|_q \mathbb{E}[Z]} + 2q \|M\|_q\}.$$

Proof.  Define

$$Z_i = \mathbb{E}\left[\sup_{f \in \mathcal{F}}\left|\sum_{j \neq i} \varepsilon_j f(X_j)\right|\,\bigg|\,X_1^n\right].$$

The monotonicity of conditional Rademacher averages with respect to the sequence of summands is well known, as it was at the core of the early concentration inequalities used in the theory of probability in Banach spaces (see [29]). Thus, for all $i$, $Z - Z_i \geq 0$ and

$$\sum_i (Z - Z_i) \leq Z.$$

Thus, we have

$$V \leq ZM \quad \text{and} \quad Z - Z_i \leq M.$$

The result now follows by Corollary 3, noticing that $M = W$.  □

**9. Moment inequalities for Rademacher chaos.**  Throughout this section, $X_1, X_2, \ldots, X_n$ denote independent Rademacher random variables. Let $\mathcal{I}_{n,d}$ be the family of subsets of $\{1, \ldots, n\}$ of size $d$ $(d < n)$. Let $\mathcal{T}$ denote a set of vectors indexed by $\mathcal{I}_{n,d}$. $\mathcal{T}$ is assumed to be a compact subset of $\mathbb{R}^{\binom{n}{d}}$.

In this section we investigate suprema of Rademacher chaos indexed by $\mathcal{T}$ of the form

$$Z = \sup_{t \in \mathcal{T}}\left|\sum_{I \in \mathcal{I}_{n,d}}\left(\prod_{i \in I} X_i\right) t_I\right|.$$

For each $1 \leq k \leq d$, let $W_k$ be defined as

$$W_k = \sup_{t \in \mathcal{T}} \sup_{\alpha^{(1)}, \ldots, \alpha^{(k)} \in \mathbb{R}^n : \|\alpha^{(h)}\|_2 \leq 1, h \leq k}$$

$$\left|\sum_{J \in \mathcal{I}_{n,d-k}}\left(\prod_{j \in J} X_j\right)\left(\sum_{i_1, \ldots, i_k : \{i_1, \ldots, i_k\} \cup J \in \mathcal{I}_{n,d}}\left(\prod_{h=1}^{k} \alpha_{i_h}^{(h)}\right) t_{\{i_1, \ldots, i_k\} \cup J}\right)\right|.$$

(Note that $W_d$ is just a constant, and does not depend on the value of the $X_i$'s.) The main result of this section is the following.

THEOREM 14.  *Let $Z$ denote the supremum of a Rademacher chaos of order $d$ and let $W_1, \ldots, W_d$ be defined as above. Then for all reals $q \leq 2$,*

$$\|(Z - \mathbb{E}[Z])_+\|_q \leq \sum_{j=1}^{d-1}(4\kappa q)^{j/2}\mathbb{E}[W_j] + (4\kappa)^{(d-1)/2}\sqrt{2K}\,q^{d/2}W_d$$

$$\leq \sum_{j=1}^{d}(4\kappa q)^{j/2}\mathbb{E}[W_j].$$

Before proving the theorem, we show how it can be used to obtain exponential bounds for the upper tail probabilities. In the special case of $d = 2$ we recover an inequality proved by Talagrand [40].

COROLLARY 4.    *For all $t \geq 0$,*

$$\mathbb{P}\{Z \geq \mathbb{E}[Z] + t\} \leq 2\exp\left(-\frac{\log(2)}{4\kappa}\bigwedge_{j=1}^{d}\left(\frac{t}{2d\mathbb{E}[W_j]}\right)^{2/j}\right).$$

PROOF.    By Theorem 14, for any $q$,

$$\mathbb{P}\{Z \geq \mathbb{E}[Z] + t\} \leq \frac{\mathbb{E}[(Z - \mathbb{E}[Z])_+]^q}{t^q}$$

$$\leq \left(\frac{\sum_{j=1}^{d}(4\kappa q)^{j/2}\mathbb{E}[W_j]}{t}\right)^q.$$

The right-hand side is at most $2^{-q}$ if for all $j = 1, \ldots, d$, $(4\kappa q)^{j/2}\mathbb{E}[W_j] \leq t/(2d)$. Solving this for $q$ yields the desired tail bound. $\square$

PROOF OF THEOREM 14.    The proof is based on a simple repeated application of Theorem 2. First note that the case $d = 1$ follows from Theorem 10. Assume that $d > 1$. By Theorem 2,

$$\|(Z - \mathbb{E}[Z])_+\|_q \leq \sqrt{2\kappa q}\|\sqrt{V^+}\|_q.$$

Now straightforward calculation shows that

$$\sqrt{V^+} \leq \sqrt{2}W_1$$

and therefore

$$\|(Z - \mathbb{E}[Z])_+\|_q \leq \sqrt{2\kappa q}\sqrt{2}(\mathbb{E}[W_1] + \|(W_1 - \mathbb{E}[W_1])_+\|_q).$$

To bound the second term on the right-hand side, we use, once again, Theorem 2. Denoting the random variable $V^+$ corresponding to $W_1$ by $V_1^+$, we see that $V_1^+ \leq 2W_2^2$, so we get

$$\|(W_1 - \mathbb{E}[W_1])_+\|_q \leq \sqrt{2\kappa q}\sqrt{2}(\mathbb{E}[W_2] + \|(W_2 - \mathbb{E}[W_2])_+\|_q).$$

We repeat the same argument. For $k = 1, \ldots, d-1$, let $V_k^+$ denote the variable $V^+$ corresponding to $W_k$. Then

$$V_k^+ \leq 2\sup_{t \in \mathcal{T}}\sup_{\alpha^{(1)}, \ldots, \alpha^{(k)} \in \mathbb{R}^n : \|\alpha^{(h)}\|_2 \leq 1, 1 \leq h \leq k}$$

$$\sum_i\left(\sum_{J \in I_{n,d-k} : i \in J}\left(\prod_{j \in J \setminus \{i\}} X_j\right)\left(\sum_{i_1, \ldots, i_k : \{i_1, \ldots, i_k\} \cup J \in \mathcal{I}_{n,d}}\prod_{h=1}^{k}\alpha_{i_h}^{(h)}t_{\{i_1, \ldots, i_k\} \cup J}\right)\right)^2$$

$$= 2W_{k+1}^2.$$

Thus, using Theorem 2 for each $W_k$, $k \leq d-1$, we obtain the desired inequality. $\qquad\square$

REMARK. Here we consider the special case $d = 2$. Let $\mathcal{T}$ denote a set of symmetric matrices with zero diagonal entries. The set $\mathcal{T}$ defines a Rademacher chaos of order 2 by

$$Z = 2 \sup_{t \in \mathcal{T}} \left| \sum_{i \neq j} X_i X_j t_{\{i,j\}} \right|.$$

Let $Y$ be defined as

$$Y = \sup_{t \in \mathcal{T}} \sup_{\alpha : \|\alpha\|_2 \leq 1} \sum_{i=1}^{n} X_i \sum_{j \neq i} \alpha_j t_{i,j}$$

and let $B$ denote the supremum of the $L_2$ operator norms of matrices $t \in \mathcal{T}$. Theorem 14 implies the following moment bound for $q \geq 2$:

$$\|(Z - \mathbb{E}[Z])_+\|_q \leq 4\sqrt{\kappa q}\, \mathbb{E}[Y] + 4\sqrt{2}\sqrt{\kappa K} q B.$$

By Corollary 4, this moment bound implies the following exponential upper tail bound for $Z$:

$$\mathbb{P}\{Z \geq \mathbb{E}[Z] + t\} \leq 2 \exp\left( -\log(2) \frac{t^2}{64\kappa \mathbb{E}[Y]^2} \wedge \frac{t}{16\sqrt{2}\sqrt{\kappa K} B} \right).$$

This is equivalent to Theorem 17 in [7] and matches the upper tail bound stated in Theorem 1.2 in [40]. Note, however, that with the methods of this paper we do not recover the corresponding lower tail inequality given by Talagrand.

We finish this section by pointing out that a version of Bonami's inequality [5] for Rademacher chaos of order $d$ may also be recovered using Theorem 14.

COROLLARY 5. *Let $Z$ be a supremum of Rademacher chaos of order $d$. Then*

$$(9.1) \qquad \|Z\|_q \leq \frac{\sqrt{4\kappa q d}^{\,d+1} - 1}{\sqrt{4\kappa q d} - 1} \|Z\|_2.$$

Note that Bonami's inequality states $\|Z\|_q \leq (q-1)^{d/2}\|Z\|_2$ so that the bound obtained by Theorem 14 has an extra factor of the order of $d^{d/2}$ in the constant. This loss in the constant seems to be an inevitable artifact of the tensorization at the basis of our arguments. On the other hand, the proof based on Theorem 2 is remarkably simple.

SKETCH OF PROOF OF COROLLARY 5. By Theorem 14, it suffices to check that for all $j, 1 \leq j \leq d$,

$$\mathbb{E}[W_j] \leq d^{j/2}\|Z\|_2.$$

Letting $W_0 = Z$, the property obviously holds for $j = 0$. Thus, it is enough to prove that for any $k > 1$,

$$\mathbb{E}[W_k] \leq \|W_k\|_2 \leq \sqrt{d}\|W_{k-1}\|_2.$$

To this end, it suffices to notice that, on the one hand,

$$
\begin{aligned}
\|W_k\|_2^2 = \mathbb{E}\Bigg[ & \sup_{t \in \mathcal{T}} \sup_{\alpha^{(1)},\ldots,\alpha^{(k-1)}\,:\,\|\alpha^{(h)}\|_2 \leq 1, h < k} \\
& \sum_{J,J' \in \mathcal{I}_{n,d-(k-1)}} |J \cap J'| \Bigg(\prod_{j \in J} X_j\Bigg)\Bigg(\prod_{j \in J'} X_j\Bigg) \\
& \Bigg( \sum_{i_1,\ldots,i_k\,:\,\{i_1,\ldots,i_{k-1}\}\cup J \in \mathcal{I}_{n,d}} \Bigg(\prod_{h=1}^{k-1}\alpha_{i_h}^{(h)}\Bigg) t_{\{i_1,\ldots,i_{k-1}\}\cup J}\Bigg) \\
& \Bigg( \sum_{i_1,\ldots,i_k\,:\,\{i_1,\ldots,i_{k-1}\}\cup J' \in \mathcal{I}_{n,d}} \Bigg(\prod_{h=1}^{k-1}\alpha_{i_h}^{(h)}\Bigg) t_{\{i_1,\ldots,i_{k-1}\}\cup J'}\Bigg) \Bigg]
\end{aligned}
$$

(the cumbersome but pedestrian proof of this identity is omitted), and on the other hand,

$$
\begin{aligned}
\|W_{k-1}\|_2^2 = \mathbb{E}\Bigg[ & \sup_{t \in \mathcal{T}} \sup_{\alpha^{(1)},\ldots,\alpha^{(k-1)}\,:\,\|\alpha^{(h)}\|_2 \leq 1, h < k} \\
& \sum_{J,J' \in \mathcal{I}_{n,d-(k-1)}} \Bigg(\prod_{j \in J} X_j\Bigg)\Bigg(\prod_{j \in J'} X_j\Bigg) \\
& \Bigg( \sum_{i_1,\ldots,i_k\,:\,\{i_1,\ldots,i_{k-1}\}\cup J \in \mathcal{I}_{n,d}} \Bigg(\prod_{h=1}^{k-1}\alpha_{i_h}^{(h)}\Bigg) t_{\{i_1,\ldots,i_{k-1}\}\cup J}\Bigg) \\
& \Bigg( \sum_{i_1,\ldots,i_k\,:\,\{i_1,\ldots,i_{k-1}\}\cup J' \in \mathcal{I}_{n,d}} \Bigg(\prod_{h=1}^{k-1}\alpha_{i_h}^{(h)}\Bigg) t_{\{i_1,\ldots,i_{k-1}\}\cup J'}\Bigg) \Bigg].
\end{aligned}
$$

Noticing that the contraction principle for Rademacher sums (see [29], Theorem 4.4) extends to Rademacher chaos in a straightforward way, and using the fact that $|J \cap J'| \leq d$, we get the desired result. $\square$

**10. Boolean polynomials.** The suprema of Rademacher chaos discussed in the previous section may be considered as special cases of suprema of $U$-processes. In this section we consider another family of $U$-processes, defined by bounded-degree polynomials of independent $\{0,1\}$-valued random variables. An important special case is the thoroughly studied problem of the number of occurrences of small subgraphs in a random graph.

In this section $X_1,\ldots,X_n$ denote independent $\{0,1\}$-valued random variables. Just like in the previous section, $\mathcal{I}_{n,d}$ denotes the set of subsets of size $d$ of $\{1,\ldots,n\}$ and $\mathcal{T}$ denotes a compact set of nonnegative vectors from $\mathbb{R}^{\binom{n}{d}}$. Note that in many applications of interest, for example, in subgraph-counting problems, $\mathcal{T}$ is reduced to a single vector.

The random variable $Z$ is defined as

$$Z = \sup_{t \in \mathcal{T}} \sum_{I \in \mathcal{I}_{n,d}} \left( \prod_{i \in I} X_i \right) t_I.$$

For the case $d = 1$, moment bounds for $Z$ follow from Theorem 11. For $k = 0,1,\ldots,d-1$, let $M_k$ be defined as

$$\max_{J \in \mathcal{I}_{n,d-k}} \sup_{t \in \mathcal{T}} \sum_{I \in \mathcal{I}_{n,d} : J \subseteq I} \left( \prod_{j \in I \setminus J} X_j \right) t_I.$$

Note that all $M_k$ are again suprema of nonnegative boolean polynomials, but the degree of $M_k$ is $k \le d$.

Lower tails for Boolean polynomials are by now well understood thanks to the Janson–Suen inequalities [17, 38]. On the other hand, upper tails for such simple polynomials are notoriously more difficult; see [20] for a survey. We obtain the following general result.

THEOREM 15. *Let $Z$ and $M_k$ be defined as above. For all reals $q \ge 2$,*

$$\|(Z - \mathbb{E}[Z])_+\|_q$$
$$\le 2 \sum_{j=1}^{d} \left\{ (\kappa q)^{j/2} \sqrt{\frac{d!}{(d-j)!}} \sqrt{\mathbb{E}[Z] \mathbb{E}[M_{d-j}]} + (\kappa q)^j \frac{d!}{(d-j)!} \mathbb{E}[M_{d-j}] \right\}.$$

PROOF. The proof is based on a repeated application of Corollary 3, very much in the spirit of the proof of Theorem 14. For each $i \in \{1,\ldots,n\}$, define

$$Z_i = \sup_{t \in \mathcal{T}} \sum_{I \in \mathcal{I}_{n,d} : i \notin I} t_I \left( \prod_{j \in I} X_j \right).$$

The nonnegativity assumption for the vectors $t \in \mathcal{T}$ implies that $Z_i \le Z$. Moreover,

$$Z - Z_i \le M_{d-1}$$

and
$$\sum_i (Z - Z_i) \le dZ.$$

Thus, $V \le dM_{d-1}Z$. Hence, by Corollary 3,
$$\|(Z - \mathbb{E}[Z])_+\|_q \le \sqrt{2\kappa qd\|M_{d-1}\|_q\mathbb{E}[Z]} + \kappa dq\|M_{d-1}\|_q.$$

We may repeat the same reasoning to each $M_k$, $k = d-1, \ldots, 1$, to obtain
$$\|M_k\|_q \le 2\left\{\mathbb{E}[M_k] + \frac{\kappa kq}{2}\|M_{k-1}\|_q\right\}.$$

By induction on $k$, we get
$$\|M_{d-1}\|_q \le 2\left\{\sum_0^{k=d-1} (\kappa q)^{d-1-k}\frac{(d-1)!}{k!}\mathbb{E}[M_k]\right\},$$

which completes the proof.  □

REMARK.  Just like in the case of Rademacher chaos, we may easily derive an exponential upper tail bound. By a similar argument to Corollary 4, we get
$$\mathbb{P}\{Z \ge \mathbb{E}[Z] + t\}$$
$$\le \exp\left(-\frac{\log 2}{d\kappa}\bigwedge_{1 \le j \le d}\left(\frac{t}{4d\sqrt{\mathbb{E}[Z]\mathbb{E}[M_j]}}\right)^{2/j}\wedge\left(\frac{t}{4d\mathbb{E}[M_{d-j}]}\right)^{1/j}\right).$$

REMARK.  If $d = 1$, Theorem 15 provides Bernstein-like bounds, and is, in a sense, optimal. For higher order, naive applications of Theorem 15 may not lead to optimal results. Moment growth may actually depend on the special structure of $\mathcal{T}$. Consider the prototypical triangle counting problem (see [18] for a general introduction to the subgraph counting problem).

In the $\mathcal{G}(n, p)$ model, a random graph of $n$ vertices is generated in the following way: for each pair $\{u, v\}$ of vertices, an edge is inserted between $u$ and $v$ with probability $p$. Edge insertions are independent. Let $X_{u,v}$ denote the Bernoulli random variable that is equal to 1 if and only if there is an edge between $u$ and $v$. Three vertices $u, v$ and $w$ form a triangle if $X_{u,v} = X_{v,w} = X_{w,z} = 1$. In the triangle counting problem, we are interested in the number of triangles
$$Z = \sum_{\{u,v,w\}\in\mathcal{I}_{n,3}} X_{u,v}X_{v,w}X_{u,w}.$$

Note that for this particular problem,
$$M_1 = \sup_{\{u,v\}\in\mathcal{I}_{n,2}}\sum_{w\,:\,w\notin\{u,v\}} X_{u,w}X_{v,w}.$$

$M_1$ is thus the maximum of $\binom{n}{2}$ (correlated) binomial random variables with parameters $n - 2$ and $p^2$. Applying Corollary 1 (with $A = 2$), we get for $q \geq 1$,

$$\|M_1\|_q \leq n \wedge (\mathbb{E}[M_1] + q - 1).$$

Simple computations reveal that

$$\mathbb{E}[M_1] \leq 2(\log n + np^2).$$

Applying Corollary 3 to $Z$, we finally get

$$\|(Z - \mathbb{E}[Z])_+\|_q \leq \sqrt{6\kappa q \mathbb{E}[M_1]\mathbb{E}[Z]}$$
$$+ q(\sqrt{6\kappa \mathbb{E}[Z]} + 3\kappa(n \wedge (\mathbb{E}[M_1] + 3(q-1)))),$$

which represents an improvement over what we would get from Theorem 15, and provides exponential bounds with the same flavor as those announced in [7]. However, the above inequality is still not optimal. In the following discussion we focus on upper bounds on $\mathbb{P}\{Z \geq 2\mathbb{E}[Z]\}$ when $p > \log n/n$. The inequality above, taking $q = \lfloor \frac{n^2 p^3}{288\kappa} \rfloor$ or $q = \lfloor \frac{\sqrt{\mathbb{E}[Z]}}{12\sqrt{\kappa}} \rfloor$, implies that for sufficiently large $n$,

$$\mathbb{P}\{Z \geq 2\mathbb{E}[Z]\} \leq \exp\left(-\log \frac{4}{3}\frac{n^2 p^3}{144\kappa} \vee \log 2 \frac{\sqrt{\mathbb{E}[Z]}}{12\sqrt{\kappa}}\right).$$

Recent work by Kim and Vu [21] show that better, and in a sense optimal, upper bounds can be obtained with some more work; see also [19] for related recent results. Kim and Vu use two ingredients in their analysis. In a first step, they tailor Bernstein's inequality for adequately stopped martingales to the triangle counting problem. This is not enough since it provides bounds comparable to the above inequality. In the martingale setting, this apparent methodological weakness is due to the fact that the quadratic variation process $\langle Z \rangle$ associated with $Z$ may suffer from giant jumps [larger than $\Theta(n^2 p^2)$] with a probability that is larger than $\exp(-\Theta(n^2 p^2))$. In the setting advocated here, huge jumps in the quadratic variation process are reflected in huge values for $M_1$. [In fact, the probability that $M_1 \geq np$ is larger than the probability that a single binomial random variable with parameters $n$ and $p^2$ is larger than $np$ which is larger than $\exp(-\Theta(np))$.] In order to get the right upper bound, Kim and Vu suggest a partitioning device. An edge $(u, v)$ is said to be good if it belongs to less than $np$ triangles. A triangle is good if its three edges are good. Let $Z^g$ and $Z^b$ denote the number of good and bad triangles. In order to bound the probability that $Z$ is larger than $2\mathbb{E}[Z]$, it suffices to bound the probability that $Z^g \geq 3/2\mathbb{E}[Z]$ and that $Z^b > \mathbb{E}[Z]/2$. Convenient moment bounds for $Z^g$ can be obtained easily using the main theorems of this paper. Indeed $Z^g/np$ satisfies the conditions of Corollary 1 with $A = 3$. Hence,

$$\|(Z^g - \mathbb{E}[Z^g])_+\|_q \leq \sqrt{\kappa}\left[\sqrt{3qnp\mathbb{E}[Z^g]} + \frac{3qnp}{\sqrt{2}}\right].$$

This moment bound implies that

$$\mathbb{P}\left\{Z^g \geq \frac{3}{2}\mathbb{E}[Z]\right\} \leq \exp\left(-\log\frac{4}{3}\frac{n^2p^2}{144}\right).$$

We refer the reader to ([21], Section 4.2) for a proof that $\mathbb{P}\{Z^b \geq \mathbb{E}[Z]/2\}$ is upper bounded by $\exp(-\Theta(n^2p^2))$.

The message of this remark is that (infamous) upper tail bounds concerning multilinear Boolean polynomials that can be obtained using Bernstein inequalities for stopped martingales can be recovered using the moment inequalities stated in the present paper. However, to obtain optimal bounds, subtle ad hoc reasoning still cannot be avoided.

## APPENDIX

**A.1. Modified $\phi$-Sobolev inequalities.** Recall the notation used in Section 3. As pointed out in [25], provided that $\phi''$ is strictly positive, the condition $1/\phi''$ concave is necessary for the tensorization property to hold. Here we point out the stronger property that the concavity of $1/\phi''$ is a necessary condition for the $\phi$-entropy $H_\phi$ to be convex on the set $\mathbb{L}_\infty^+(\Omega, \mathcal{A}, \mathbb{P})$ of bounded and nonnegative random variables.

PROPOSITION A.1. *Let $\phi$ be a strictly convex function on $\mathbb{R}_+$ which is twice differentiable on $\mathbb{R}_+^*$. Let $(\Omega, \mathcal{A}, \mathbb{P})$ be a rich enough probability space in the sense that $\mathbb{P}$ maps $\mathcal{A}$ onto $[0, 1]$. If $H_\phi$ is convex on $\mathbb{L}_\infty^+(\Omega, \mathcal{A}, \mathbb{P})$, then $\phi''(x) > 0$ for every $x > 0$ and $1/\phi''$ is concave on $\mathbb{R}_+^*$.*

PROOF. Let $\theta \in [0, 1]$ and $x, x', y, y'$ be positive real numbers. Under the assumption on the probability space we can define a pair of random variables $(X, Y)$ to be $(x, y)$ with probability $\theta$ and $(x', y')$ with probability $(1 - \theta)$. Then the convexity of $H_\phi$ means that

$$H_\phi(\lambda X + (1 - \lambda)Y) \leq \lambda H_\phi(X) + (1 - \lambda)H_\phi(Y)$$

for every $\lambda \in (0, 1)$. Defining, for every $(u, v) \in \mathbb{R}_+^* \times \mathbb{R}_+^*$,

$$F_\lambda(u, v) = -\phi(\lambda u + (1 - \lambda)v) + \lambda\phi(u) + (1 - \lambda)\phi(v),$$

the inequality is equivalent to

$$F_\lambda(\theta(x, y) + (1 - \theta)(x', y')) \leq \theta F_\lambda(x, y) + (1 - \theta)F_\lambda(x', y').$$

Hence, $F_\lambda$ is convex on $\mathbb{R}_+^* \times \mathbb{R}_+^*$. This implies, in particular, that the determinant of the Hessian matrix of $F_\lambda$ is nonnegative at each point $(x, y)$. Thus, setting $x_\lambda = \lambda x + (1 - \lambda)y$,

$$[\phi''(x) - \lambda \phi''(x_\lambda)][\phi''(y) - (1 - \lambda)\phi''(x_\lambda)] \geq \lambda(1 - \lambda)[\phi''(x_\lambda)]^2,$$

which means that

$$(A.1) \qquad \phi''(x)\phi''(y) \geq \lambda \phi''(y)\phi''(x_\lambda) + (1 - \lambda)\phi''(x)\phi''(x_\lambda).$$

If $\phi''(x) = 0$ for some point $x$, we see that either $\phi''(y) = 0$ for every $y$, which is impossible because $\phi$ is assumed to be strictly convex, or there exists some $y$ such that $\phi''(y) > 0$ and then $\phi''$ is identically equal to 0 on the nonempty open interval with endpoints $x$ and $y$, which also leads to a contradiction with the assumption that $\phi$ is strictly convex. Hence $\phi''$ is strictly positive at each point of $\mathbb{R}_+^*$ and (A.1) leads to

$$\frac{1}{\phi''(\lambda x + (1 - \lambda)y)} \geq \frac{\lambda}{\phi''(x)} + \frac{1 - \lambda}{\phi''(y)},$$

which means that $1/\phi''$ is concave. □

PROOF OF PROPOSITION 2. Without loss of generality we may assume that $\phi(0) = 0$. If $\phi$ is strictly convex,

$$\frac{1}{\phi''((1 - \lambda)u + \lambda x)}$$

$$\geq \frac{1 - \lambda}{\phi''(u)} + \frac{\lambda}{\phi''(x)} \qquad \text{(by concavity of } 1/\phi'')$$

$$\geq \frac{\lambda}{\phi''(x)} \qquad \qquad \text{(by positivity of } \phi'', \text{ i.e., strict convexity of } \phi).$$

In any case, the concavity of $1/\phi''$ implies that for every $\lambda \in (0, 1)$ and every positive $x$ and $u$,

$$\lambda \phi''((1 - \lambda)u + \lambda x) \leq \phi''(x),$$

which implies that for every positive $t$,

$$\lambda \phi''(t + \lambda x) \leq \phi''(x).$$

Letting $\lambda$ tend to 1, we derive from the above inequality that $\phi''$ is nonincreasing, that is, $\phi'$ is concave. Setting $\psi(x) = \phi(x)/x$, one has

$$x^3 \psi''(x) = x^2 \phi''(x) - 2x \phi'(x) + 2\phi(x) = f(x).$$

The convexity of $\phi$ and its continuity at 0 imply that $x\phi'(x)$ tends to 0 as $x$ goes to 0. Also, the concavity of $\phi'$ implies that

$$x^2 \phi''(x) \leq 2x(\phi'(x) - \phi'(x/2)),$$

so $x^2 \phi''(x)$ tends to 0 as $x \to 0$ and therefore $f(x) \to 0$ as $x \to 0$. Denoting (abusively) by $\phi^{(3)}$ the right derivative of $\phi''$ (which is well defined since $1/\phi''$ is concave) and by $f'$ the right derivative of $f$, we have $f'(x) = x^2 \phi^{(3)}(x)$. Then $f'(x)$ is nonpositive because $\phi''$ is nonincreasing. Thus, $f$ is nonincreasing. Since $f$ tends to 0 at 0, this means that $f$ is a nonpositive function and the same property holds for the function $\psi''$, which completes the proof of the concavity of $\psi$.   $\square$

Proof of Lemma 2.  Without loss of generality we assume that $\phi(0) = 0$. The convexity of $\phi$ implies that for every positive $u$,

$$-\phi(\mathbb{E}[Z]) \leq -\phi(u) - (\mathbb{E}[Z] - u)\phi'(u),$$

and therefore

$$H_\phi(Z) \leq \mathbb{E}[\phi(Z) - \phi(u) - (Z - u)\phi'(u)].$$

Since the latter inequality becomes an equality when $u = m$, the variational formula (3.7) is proven. Since $Z'$ is an independent copy of $Z$, we derive from (3.7) that

$$H_\phi(Z) \leq \mathbb{E}[\phi(Z) - \phi(Z') - (Z - Z')\phi'(Z')]$$
$$\leq -\mathbb{E}[(Z - Z')\phi'(Z')]$$

and by symmetry

$$2H_\phi(Z) \leq -\mathbb{E}[(Z' - Z)\phi'(Z)] - \mathbb{E}[(Z - Z')\phi'(Z')],$$

which leads to (3.8). To prove (3.9), we simply note that

$$\tfrac{1}{2}\mathbb{E}[(Z - Z')(\psi(Z) - \psi(Z'))] - H_\phi(Z) = -\mathbb{E}[Z]\mathbb{E}[\psi(Z)] + \phi(\mathbb{E}[Z]).$$

But the concavity of $\psi$ implies that $\mathbb{E}[\psi(Z)] \leq \psi(\mathbb{E}[Z]) = \phi(\mathbb{E}[Z])/\mathbb{E}[Z]$ and we derive from the preceding identity that (3.9) holds.   $\square$

Proof of Theorem 6.  Fix first $\tilde{y} \leq x \leq y$. Under the assumption that $g = \phi' \circ f$ is convex,

(A.2)
$$\phi(f(y)) - \phi(f(x)) - (f(y) - f(x))\phi'(f(x))$$
$$\leq \tfrac{1}{2}(y - x)^2 f'^2(\tilde{y})\phi''(f(\tilde{y})).$$

Indeed, denoting by $h$ the function

$$h(t) = \phi(f(y)) - \phi(f(t)) - (f(y) - f(t))g(t),$$

we have

$$h'(t) = -g'(t)(f(y) - f(t)).$$

But for every $t \leq y$, the monotonicity and convexity assumptions on $f$ and $g$ yield

$$0 \leq -g'(t) \leq -g'(\tilde{y}) \quad \text{and} \quad 0 \leq -(f(y) - f(t)) \leq -(y - t)f'(\tilde{y}),$$

hence

$$-h'(t) \leq (y - t)f'(\tilde{y})g'(\tilde{y}).$$

Integrating this inequality with respect to $t$ on $[x, y]$ leads to (A.2). Under the assumption that $\psi \circ f$ is convex, we notice that

$$0 \leq -(f(y) - f(x)) \leq -(y - x)f'(\tilde{y})$$

and

$$0 \leq -(\psi(f(y)) - \psi(f(x))) \leq -(y - x)f'(\tilde{y})\psi'(f(\tilde{y})),$$

which implies

(A.3)    $$(f(y) - f(x))(\psi(f(y)) - \psi(f(x))) \leq (x - y)^2 f'^2(y)\psi'(f(y)).$$

The tensorization inequality combined with (3.7) and (A.2) leads to

$$H_\phi(f(Z)) \leq \frac{1}{2} \sum_{i=1}^{n} \mathbb{E}[(Z - Z_i)^2 f'^2(\tilde{Z}) \phi''(f(\tilde{Z}))]$$

and therefore to the first inequality of Theorem 6, while we derive from the tensorization inequality (3.9) and (A.3) that

$$H_\phi(f(Z)) \leq \sum_{i=1}^{n} \mathbb{E}[(Z - Z_i')_+^2 f'^2(\tilde{Z})\psi'(f(\tilde{Z}))],$$

which means that the second inequality of Theorem 6 indeed holds.

In order to prove the third inequality, we simply define $\tilde{f}(x) = f(-x)$ and $\tilde{Z} = -Z$. Then $\tilde{f}$ is nondecreasing and convex and we can use the first inequality of Theorem 5 to bound $H_\phi(f(Z)) = H_\phi(\tilde{f}(\tilde{Z}))$, which gives

$$H_\phi(\tilde{f}(\tilde{Z})) \leq \sum_{i=1}^{n} \mathbb{E}[(\tilde{Z} - \tilde{Z}_i')_+^2 \tilde{f}'^2(\tilde{Z})\psi'(\tilde{f}(\tilde{Z}))]$$

$$\leq \sum_{i=1}^{n} \mathbb{E}[(Z - Z_i')_-^2 f'^2(Z)\psi'(f(Z))],$$

completing the proof of the result.   $\square$

**A.2. Proof of Lemma 6.** By Stirling's formula,

$$k! = \left(\frac{k}{e}\right)^k \sqrt{2\pi k}\, e^{2\beta_k},$$

where $\beta_k$ is positive and decreases to 0 as $k \to \infty$. Using the above formula with $k = q-2$, $k = q-1$ and $k = q$ leads to

$$x_q \le e^{\beta_q - \beta_{q-1} - 1/2}\left(\frac{q-1}{q}\right)^{1/4} + \frac{1}{K} e^{\beta_q - \beta_{q-2} - 1}\left(\frac{(q-1)^2}{q(q-2)}\right)^{1/4}.$$

By the monotonicity of Stirling's correction, we have $\beta_q \le \beta_{q-1} \le \beta_{q-2}$, and the preceding inequality becomes

$$x_q \le e^{-1/2}\left(\frac{q-1}{q}\right)^{1/4} + \frac{1}{K} e^{-1}\left(\frac{q-1}{(q-2)q}\right)^{1/4}.$$

Our aim is to prove that $x_q \le 1$. Let

$$a_q = e^{-1/2}(q-1)^{1/4} q^{-1/4},$$

$$u_q = \frac{e^{-1}(q-1)^{1/2}(q-2)^{-1/4} q^{-1/4}}{1 - a_q}.$$

Then

$$x_q \le a_q + \frac{1}{K} u_q (1 - a_q)$$

and since $u_q \to K$ as $q \to \infty$, in order to show that $x_q \le 1$, it is enough to prove that $u_q \le u_{q+1}$ for every $q \ge 4$. Let $\theta = 1/q$. Then $u_q \le u_{q+1}$ is equivalent to $g(\theta) \ge 0$, where

$$\begin{aligned}
\text{(A.4)} \quad g(\theta) = {} & (1 - 2\theta)^{1/4}(1 - e^{-1/2}(1-\theta)^{1/4}) + e^{-1/2}(1-\theta)^{3/4} \\
& - (1-\theta)^{1/2}(1-\theta^2)^{1/4}.
\end{aligned}$$

Now, $t \to t^{-2}((1-t)^{1/2} - (1-2t)^{1/4})$ is easily seen to be increasing on $(0, 1/2)$ (just notice that its power series expansion has nonnegative coefficients), so since $\theta \le 1/4$, setting $\gamma = 16(\sqrt{3/4} - 2^{-1/4})$, one has $(1-2\theta)^{1/4} \ge (1-\theta)^{1/2} - \gamma\theta^2$. Plugging this inequality in (A.4) yields

$$g(\theta) \ge (1-\theta)^{1/2}(1 - (1-\theta^2)^{1/4}) - \gamma\theta^2(1 - e^{-1/2}(1-\theta)^{1/4})$$

and therefore, using again that $\theta \le 1/4$,

$$g(\theta) \ge (\tfrac{3}{4})^{1/2}(1 - (1-\theta^2)^{1/4}) - \gamma\theta^2(1 - e^{-1/2}(\tfrac{3}{4})^{1/4}).$$

Finally, note that $1 - (1-\theta^2)^{1/4} \ge \theta^2/4$ which implies

$$\theta^{-2} g(\theta) \ge \tfrac{1}{4}(\tfrac{3}{4})^{1/2} - \gamma(1 - e^{-1/2}(\tfrac{3}{4})^{1/4})$$

and since one can check numerically that the right-hand side of this inequality is positive (more precisely, it is larger than 0.041), we derive that the sequence $(u_q)$ is increasing and is therefore smaller than its limit $K$ and the result follows.

**A.3. Proof of Lemma 8.** The statement follows from a version of Hoffmann–Jørgensen's inequality. In particular, we use inequality (1.2.5s) on page 10 in [13] with $p = 1$ and $t = 0$. Then we obtain

$$\mathbb{E}\left[\sup_f \left|\sum_i \varepsilon_i f^2(X_i) \mathbb{1}_{\sup_f |f(X_i)| > t_0}\right|\right]$$
$$\leq \left(\frac{\mathbb{E}[M^2]^{1/2}}{1 - (4\mathbb{P}[\sup_f |\sum_i \varepsilon_i f^2(X_i)| > t_0])^{1/2}}\right)^2.$$

The right-hand side may be bounded further by observing that, by Markov's inequality,

$$\mathbb{P}\left[\sup_f \left|\sum_i \varepsilon_i f^2(X_i) \mathbb{1}_{\sup_f |f(X_i)| > t_0}\right| > 0\right] = \left[\sup_{f,i} |f(X_i)| > t_0\right]$$
$$\leq \frac{\mathbb{E}[M^2]}{t_0^2} = \frac{1}{\lambda}.$$

**Acknowledgment.** The authors gratefully acknowledge the care and patience of an anonymous referee.

S. BOUCHERON
LRI, UMR 8623 CNRS
CNRS-UNIVERSITÉ PARIS-SUD
BÂTIMENT 490
91405 ORSAY-CEDEX
FRANCE
E-MAIL: stephane.boucheron@lri.fr

O. BOUSQUET
MAX PLANCK INSTITUTE
   FOR BIOLOGICAL CYBERNETICS
SPEMANNSTR. 38
D-72076 TÜBINGEN
GERMANY

G. LUGOSI
DEPARTMENT OF ECONOMICS
POMPEU FABRA UNIVERSITY
RAMON TRIAS FARGAS 25-27
08005 BARCELONA
SPAIN

P. MASSART
MATHÉMATIQUES
UNIVERSITÉ PARIS-SUD
BÂTIMENT 425
91405 ORSAY-CEDEX
FRANCE